\input amstex
\documentstyle{amsppt}
\mag=\magstep2 
\pagewidth{11cm}
\NoBlackBoxes
\pageheight{15cm}
\def\a{\alpha}
\def\b{\beta}
\def\e{\varepsilon}
\def\gam{\gamma}
\def\Gam{\Gamma}
\def\k{\kappa}

\def\lam{\lambda}
\def\ome{\omega}
\def\Ome{\Omega}
\def\sig{\sigma}

\def\A{\Cal A }
\def\E{\Cal E}

\def\G{G_2}

\def\Diff{\text{\rm Diff}_0}
\def\R{\Bbb R}
\def\C{\Bbb C}
\def\H{\Bbb H}
\def\Z{\Bbb Z}
\def\O{\Cal O}

\def\M{\frak M}

\def\ol{\overline}
\def\w{\wedge}
\def\({\left(}
\def\){\right)}
\def\G{G_2}

\def\neg{\negthinspace}
\def\h{\hat}
\def\hrho{\hat{\rho}}

\def\til{\tilde}
\def\wtil{\widetilde}

\def\pa{\partial}

\def\Spin{\text{Spin}(7)}
\def\s-style{\scriptstyle}
\def\ss-style{\scriptscriptstyle}
\def\trian{\triangle}

\def\arrow{\longrightarrow}

\def\CP{\dsize{\Bbb C} \bold P}

\def\W{\wedge_{\ome^0}}
\def\sw{{\ss-style{\w}}}
\document
\topmatter
\title\nofrills Moduli spaces of topological calibrations,\\
{	Calabi-Yau, 
HyperK\"ahler,\\
$\G$ and Spin$(7)$ structures}
\endtitle 
\rightheadtext{Topological calibrations}
\author	
Ryushi Goto
\endauthor
\affil		
Department of Mathematics,\\
 Graduate School of Science,\\
Osaka University,\\ 
\endaffil
\address
Toyonaka, Osaka, 560, Japan
\endaddress
\email
goto\@math.sci.osaka-u.ac.jp
\endemail
\abstract
We shall obtain unobstructed deformations
of four geometric structures: Calabi-Yau,
HyperK\"ahler, $\G$ and Spin$(7)$ structures 
in terms of closed differential forms (calibrations).
We develop a direct and unified construction of smooth moduli spaces
of these four geometric structures
and show that the local Torelli type theorem holds in a systematic way.
\endabstract
\endtopmatter

\head 
\S0. Introduction 
\endhead
There has been considerable interest recently in Riemannian manifolds 
with vanishing Ricci tensor. 
The list of holonomy group of Ricci-flat manifolds includes four interesting 
classes of the holonomy groups: SU$(n)$, Sp$(m)$, $\G$ and $\Spin$ [1]. 
The Lie group SU$(n)$ arises as the holonomy group of Calabi-Yau
manifolds and Sp$(m)$ is the holonomy group of 
HyperK\"ahler manifolds. The exceptional Lie group $\G$ and the Lie group $\Spin$
 respectively occur as the holonomy groups of 
$7$ and $8$ dimensional manifolds.
There are many interesting common
properties between these four geometries. One of the most remarkable property is
that the deformation spaces of these geometric structures are smooth,i.e.,
unobstructed. It is also  
intriguing that there are smooth moduli spaces of these geometric structures, 
in which the local Torelli type theorem hold, so that is, 
these moduli spaces are locally described in terms of cohomology groups.
Bogomorov, Tian and Todorov show that 
the deformation space (Kuranishi space) of Calabi-Yau structures is smooth by using
Kodaira-Spencer theory [2],[25],[27]. 
The moduli space of polarized Calabi-Yau manifolds is
constructed by Fujiki-Schmacher [8] from complex geometric point of view.
Joyce obtains smooth moduli spaces
of $\G$ and
$\Spin$ structures respectively [13],[14],[15]. The construction of
moduli spaces of $\G$ and Spin$(7)$ structures are different from the one of
Calabi-Yau structures since
$\G$ and $\Spin$ manifolds are real manifolds, in which we can not apply the
deformation theory of complex manifolds. Hitchin shows a significant and suggestive
construction of deformation spaces of Calabi-Yau structures on real $6$ manifolds
and $\G$ structures on $7$ manifolds [12]. It must be noted that these four
geometries are defined by certain closed differential forms on real manifolds. 
From this point of view  we shall obtain a direct and unified construction of 
smooth moduli spaces of these geometric structures. 
In the case of Calabi-Yau manifolds, we consider a real
compact $2n$ manifold with a pair consisting of a closed complex $n$ form $\Ome$
and  a symplectic form $\ome$.  We show that a certain pair $(\Ome,\ome)$ defines a
Calabi-Yau metric ( Ricci-flat K\"ahler metric) on $X$. Hence the deformation space
of Calabi-Yau metrics on $X$ arises as the deformation space of such pairs of
closed forms 
$(\Ome,\ome)$ 
(see section 4-2 for precise definition of Calabi-Yau structures ).  In section 1,
we discuss a general deformation theory of geometric structures defined by closed
differential forms. Let $V$ be a real $n$ dimensional vector space.  Then we
consider the linear action $\rho$ of $G=$GL$(V)$ on the direct sum of
skew-symmetric tensors, 
$$
\rho \: \text{GL}(V) \to \oplus_{i=1}^l\text{End}(   \w^{p_i} V^*).
$$
Let $\Phi^0_{\ss-style V}=(\phi^0_1,\phi^0_2,\cdots, \phi^0_l)$ be an element of $\oplus_{i=1}^l
\w^{p_i} V^*$. Then we have the $G$-orbit
$\O=\O_{\ss-style\Phi^0_{\ss-style V}}(V)$:
$$
\O_{\ss-style\Phi^0_{\ss-style V}}(V) = \{\,\rho_g \Phi^0_{\ss-style V}=(\rho_g \phi^0_1, \cdots 
,\rho_g\phi^0_l ) \in
\oplus_{i=1}^l\w^{p_i}  V^*\, |\, g \in
G\,\}. 
$$
Then the orbit $\O=\O_{\ss-style\Phi^0_{\ss-style V}}$ is regarded as a homogeneous space $G/H$, 
where $H$ is the isotropy group. 
If the isotropy group H is a subgroup of the orthogonal group O$(V)$ for a metric
$g_V$ on
$V$,  we call $\O$ an metrical orbit.
Let $X$ be a real $n$ dimensional compact manifold. 
Then we define a homogeneous space bundle $\A_{\ss-style\O}(X) \to X$ by 
$$
\A_{\ss-style{\O}}(X)=\bigcup_{x\in X}\O_{\ss-style\Phi^0_{\ss-style V}}(T_x X).
$$
Then we define $\E_{\ss-style{\Cal O}}(X)$ to be the set of global sections
$\Gam(X,\A_{\ss-style\O}(X))$. 
The moduli space $\M_{\ss-style{\O}}(X)$ is defined as the quotient
$$
\M_{\ss-style{\O}}(X)=\wtil{\M}_{\ss-style{\O}}(X)/\Diff(X),
$$
where 
$$
\wtil{\M}_{\ss-style{\O}}(X)=\{\, \Phi \in \E(X)\, |\, d\Phi =0 \, \}
$$ and $\Diff(X)$ denotes the identity component of 
diffeomorphisms of $X$.  
Let $\Phi^0$ be an element of $\widetilde{\M}_{_\O} (X)$. 
Then we shall obtain a deformation complex $\#_{\Phi^0}$ (see section 1): 
$$
\CD 
\Gam(E^0)@>d_0>>\Gam(E^1)@>d_1>>\Gam(E^2)@>d_2>>\cdots
\endCD
$$
If $\#_{\Phi^0}$ is an elliptic complex, the orbit $\O$ is called an elliptic orbit 
(see definition 1-1). 
As we shall
show that this complex
$\#_{\Phi^0}$ is a subcomplex of the direct sum of de Rham complex
(for simplicity we call this the de Rham complex): 
$$
\CD 
@.\Gam(E^0)@>d_0>>\Gam(E^1)@>d_1>>\Gam(E^2)@>d_2>>\cdots \\
@. @VVV @VVV @VVV @. \\
\cdots @>d>>\Gam ( \oplus_i\w^{p_i -1}) @>d>> \Gam (\oplus_i \w^{p_i}) @>d>> 
\Gam (\oplus_i \w^{p_i +1} )@>d>>\cdots
\endCD
$$
Hence we have the map $p^k$ from the cohomology group of the complex $\#_{\Phi^0}$
to the cohomology group of the de Rham complex: 
$$
p^k\:H^k (\#_{\Phi^0}) \arrow \oplus_i H^{p_i-1 -k }_{dR}(X).
$$ 
If the maps $p^1$ and $p^2$ are respectively injective for any $\Phi^0 \in
\wtil{\M}_{\ss-style{\O}}(X)$ on every compact $n$ dimensional manifold $X$, 
we call $\O$ a topological orbit.
Then we have the following theorems:
\proclaim{Theorem 1-8 } 
If an orbit $\O$ is metrical, elliptic and topological, 
then the corresponding moduli space $\M_{_\O} (X)$ is a smooth manifold. 
( In particular $\M_{_\O} (X)$ is Hausdorff .)
Further $\M_{_\O} (X)$ has canonical coordinates given by an open ball of the
cohomology group H$^1 (\#_{\Phi}) $. 
\endproclaim
Since de Rham cohomology group is invariant under the action of Diff$_0(X)$, we
have the map 
$$
P\: \M_{_\O} (X) \arrow \underset i \to \oplus H^{p_i}_{dR} (X).
$$
Then we have 
\proclaim{Theorem 1-9}
If an orbit $\O$ is metrical, elliptic and topological, then the map $P$ is locally injective. 
\endproclaim
Further under the assumption that $\O $ is metrical, elliptic and topological 
we have
\proclaim{Theorem 1-11}
Let $\wtil{\M}_{\ss-style\O}(X)$ be the set of closed elements of $\E$. 
We denote by {\rm Diff}$(X)$ the group of diffeomorphisms of $X$. There is the
action of {\rm Diff}$(X)$ on $ \wtil{\M}_{\ss-style\O}(X)$.  Then the quotient
$\wtil{\M}_{\ss-style\O}(X)/\text{\rm Diff}(X)$ is an orbifold.  
\endproclaim
In order to obtain these theorems, we study the problem of obstruction 
to deformations in our situation. 
$\E_{\ss-style{\Cal O}}(X)$ is regarded as
 a infinite dimensional homogenous space (a Hilbert manifold). Hence we have the
tangent space 
$T_{\ss-style\Phi^0}\E_{\ss-style{\Cal O}}(X)$ of $\E_{\ss-style{\Cal O}}(X)$. We
denote by $\Cal H$ the Hilbert space consisting of closed forms in $\oplus_i
\w^{p_i}$ Then the space
$\wtil{\M}_{\ss-style{\O}}(X)$ is the intersection  between the Hilbert space $\Cal
H$ and the Hilbert manifold $\E{\ss-style{\Cal O}}(X)$.  We define an infinitesimal
tangent space  of
$\wtil{\M}_\O$ by the intersection 
$\Cal H\cap T_{\ss-style\Phi^0}\E_{\ss-style{\Cal O}}(X)$.
Then we shall discuss if the infinitesimal tangent space
is regarded as the tangent space of actual deformations.
\proclaim{Definition 1-6} 
A closed element $\Phi^0 \in \E(X)$ is unobstructed if there 
exists an integral curve $\Phi_t(\a)$ in $\wtil{\M}_{\ss-style{\Cal O}}(X)$ 
for each infinitesimal tangent vector 
$\a\in \Cal H\cap T_{\ss-style\Phi^0}\E_{\ss-style{\Cal O}}(X) $
such  that 
$$
\frac d{dt}\Phi_t(\a)|_{t=0} = \a
$$
An orbit $\O$ is unobstructed if every $\Phi^0 \in \wtil{\M}_{\O}(X)$  is
unobstructed  for any compact $n$ dimensional manifold $X$
\endproclaim
We shall prove the following criterion in section 2.
\proclaim{Theorem 1-7 ( Criterion of unobstructedness)} 
We assume that an orbit $\O$ is elliptic (see definition 1-1 in section one). 
If the map $p^2\:H^2(\#_{\ss-style\Phi^0})\to \oplus_i H^{p_i+1}_{\ss-style DR}(X)$
is injective, then $\Phi^0$ is unobstructed ( see section 1 for $p^2$).
\endproclaim
At first we try to construct a deformation of calibrations as a formal
power series in $t$.
Then we encounter obstructions to deformation of calibrations.  
A primary obstruction is discussed in subsection 2-1, 
which is given by a generalization of the Nijenhuis tensor (see subsection 2-0). If
the primary obstruction vanishes, then we have the second obstruction. Successively
we have higher obstructions to deformations. Explicit description of higher
obstructions are given in subsection 2-2. In subsection 2-3, we prove our criterion
of unobstructedness (Theorem 1-5).
If the criterion holds, then all obstruction vanish simultaneously.
Hence we have a deformation of calibrations as a formal
power series in $t$. Further we prove the power series uniformly converges.
Section 3 is devoted to prove main theorems. 
Our discussion is based on 
[7] and [22].
We have a smooth family of closed forms $\Cal S_{\Phi^0}$, 
which is given by the deformation space in section 2. 
The injectivety of the map $p^1$ is essential to show that $\Cal S_{\Phi^0}$ 
gives coordinates of the moduli space $\M_{\ss-style{\O}}(X)$
in subsection 3-2. We also obtain the Hausdorff
property of the moduli space in subsection 3-3. 
We will give the proof of theorems in subsection 3-4. 
In section $4,5,6$ and $7$ we shall show that 
Calabi-Yau, HyperK\"ahler, $\G$ and Spin$(7)$ structures are metrical, elliptic 
and topological respectively. 
In section 4-1 we define an {\it SL$_n(\C)$ structure } as a certain complex
form $\Ome$, which defines the almost complex structure $I_{\Ome}$ with trivial
canonical line bundle.  Then the integrability of the almost complex structure 
$I_{\Ome}$ is given by a closeness of the complex differential form $\Ome$.
We show that 
the orbit of SL$_n(\C)$ structures is elliptic and 
satisfies the criterion. In section
4-2, we define a {\it Calabi-Yau structure} as a certain pair consisting 
SL$_n(\C)$ structure $\Ome$ and a real symplectic form $\ome$.
Then we
prove that the orbit corresponding to a Calabi-Yau structure is
elliptic and  topological. 
Hence we obtain the smooth moduli space of Calabi-Yau structures. 
A reference of Calabi-Yau manifolds is [1]. 
Our primary obstruction of SL$_n(\C)$ structures corresponds to the one of
Kodaira-Spencer theory. Then our result is regarded as 
another proof of unobstructedness  
by using calibrations. 
Our direct proof reveals a geometric meaning of unobstructed deformations. 
( we do not use Calabi-Yau's theorem to 
obtain a smooth deformation space of Calabi-Yau structures). 
It must be noted that Kawamata and Ran give algebraic proof of unobstructed
deformations.)[16],[23]. 
In section
5, we show the orbit corresponding to a HyperK\"ahler structure is also
elliptic and topological. In section 6 and 7 we discuss $\G$ and Spin $(7)$
structures respectively.

\head \S1. Moduli spaces of calibrations
\endhead
Let $V$ be a real vector space of dimension $n$. We denote by $\w^p V^*$ the
vector space  of $p$ forms on $V$. Let $\rho_p$ be the linear action of
$G=$GL$(V)$ on $\w^p V^*$. Then we have the action $\rho$ of $G$ on 
the direct sum
$\oplus_i
\w^{p_i} V^*$ by 
$$
\rho \: GL(V) \arrow \oplus_{i=1}^l\text{End}( \w^{p_i} V^*),
$$
$$\rho= (\rho_{p_1},\cdots,\rho_{p_l} ).
$$
We fix an element $\Phi^0_{\ss-style V} = ( \phi^0_1, \phi^0_2, \cdots ,\phi^0_l ) \in
\oplus_i
\w^{p_i}V^*$  and consider the $G$-orbit $\O=\Cal O_{\Phi^0_{\ss-style V}}$ through
$\Phi^0_{\ss-style V}$: 
$$
\Cal O_{\Phi^0_{\ss-style V}} = 
\{ \, \Phi_{\ss-style V} = \rho_g \Phi^0_{\ss-style V} \in \oplus_i \w^{p_i}
V^*
\, |\, g
\in G\, \}
$$
The orbit $\O_{\Phi^0_{\ss-style V}}$ can be regarded as a homogeneous space,
$$
\Cal O_{\Phi^0_{\ss-style V}} = G/ H,
$$where $H$ is the isotropy group 
$$
H = \{ \, g \in G \, | \, \rho_g \Phi^0_{\ss-style V} = \Phi^0_{\ss-style V} \, \}.
$$
We denote by $\Cal A_\O (V) =\Cal A (V)$ the orbit $\Cal O_{\Phi^0_{\ss-style
V}} = G/H$.  The tangent space $E^1 (V)=T_{\Phi^0_{\ss-style V}} \A(V)$ is
given by 
$$
E^1(V)=T_{\Phi^0_{\ss-style V}}\A(V) = 
\{\, \rho_\xi \Phi^0 \in \oplus_i \w^{p_i}V^* \, |\, \xi \in \frak g\, \},
$$
where $\rho$ denotes the differential representation of $\frak g$.
The vector space $E^1 (V)$ is the quotient space $\frak g/\frak h$.
We also define a vector space $E^0(V)$ by 
the interior product, 
$$\align
E^0(V) =& \{\, i_v \Phi^0_{\ss-style V} =(i_v\phi^0_1,\cdots,i_v\phi^0_l)\in \oplus_i
\w^{p_i-1} V^*\, |\, v
\in V
\,
\}.
\endalign
$$
$E^2(V)$ is define as a vector space spanned by the following set, 
$$
E^2(V) =\text{Span}\{\, \a\w i_v\Phi^0_{\ss-style V} \in \oplus_i \w^{p_i +1} V^* \, |\,
\a \in \w^2V^* , \, 
i_v\Phi^0 \in E^0 (V) \, \}.
$$
We also define $E^k(V)$ for $k \geq 0$ by 
$$
E^k(V) = \text{Span}\{\, \b\w i_v\Phi^0_{\ss-style V} \in \oplus_i \w^{p_i +k-1} V^* \, |\,
\b \in \w^{k} V^* , \, 
i_v\Phi^0_{\ss-style V} \in E^0 (V) \, \}.
$$
Let$\{ e_1, \cdots , e_n \} $   be a basis
of $V$ and  $\{\theta^1, \cdots , \theta^n \, \}$ the dual basis of $V^*$.
Then we see that $\rho_\xi \Phi^0_{\ss-style V}$ is written as 
$$\rho_{\xi}\Phi^0_{\ss-style V} = \sum_{i j} \xi_i ^j \theta^j \w
i_{e_i} \Phi^0_{\ss-style V},$$
where $\xi = \sum_{ij}\xi_j^i \theta^j \otimes e_i $ and $i_{e_i}$ denotes the
interior product. Hence we have the graded vector space $E(V) = \oplus_k E^k
(V)$ generated by $E^0(V)$ over $\w^*V^* $.
Then we have the complex by the exterior product of a nonzero $u \in V^*$, 
$$
\CD
E^0 (V) @>\w u>>E^1(V)@>\w u >>E^2 (V)@>\w u>>\cdots. 
\endCD
$$
\proclaim{Definition 1-1(elliptic orbits)}
An orbit $\O_{\Phi^0_{\ss-style V}}$ is an elliptic orbit if the complex 
$$
\CD
E^0 (V) @>\w u>>E^1(V)@>\w u >>E^2 (V)@>\w u>>\cdots. 
\endCD
$$
is exact for any nonzero $u \in V^*.$ 
In other words, if $\a\w u =0$ for $\a\in E^k(V)$, then there exists $\b \in
E^{k-1}(V)$ such that $\a = \b \w u$ for $k=1,2$.
\endproclaim
\demo{Remark}
If $\a\w u =0$, then
we have $\a = \b \w u$
for some $\b \in \oplus_i \w^{p_i -1}$ since the  de Rham complex is
elliptic. However 
$\b$ is not an element of $E^0(V)$ in general, 
(Note that $E^0(V)$ is a subspace of $\oplus_i \w^{p_i -1}$).
For instance, we take
$\Phi^0_{\ss-style V}$ as a real symplectic form $\ome$  on a real $2n$
dimensional vector space
$V$.  Then $E^0 =\w^1$ and $E^1 = \w^2$. Hence $\Cal O_{\Phi^0_{\ss-style V}}$
is  elliptic. However if $\Phi^0_{\ss-style V}$ is a degenerate $2$ form on $V$,
i.e., $\ome^n =0$, then  
$\Cal O_{\Phi^0_{\ss-style V}}$ is not elliptic. 
\enddemo
\proclaim{Definition 1-2(metrical orbits)}
Let $\O_{\ss-style\Phi^0_{\ss-style V}}$ be an orbit as before.
An orbit $\O_{\ss-style\Phi^0_{\ss-style V}}$ is metrical if the isotropy group
$H$ is a subgroup of  O$(V)$ with respect to a metric $g_V$ on $V$.
\endproclaim
Let $X$ be a compact real manifold of dimension $n$. 
We define $\A_{\ss-style \O}(T_x X)$ by using an identification 
$h\: T_xX \cong V$. The subspace $\A_{\ss-style \O}(T_x X) \subset 
\oplus_i\w^{p_i}T^*_x X$ is independent of a choice of an identification $h$.
Hence we define the $G/H-$bundle $\A(X)(=\A_{\ss-style{\Cal O}}(X))$ by 
$$
\A_{\ss-style\O}(X) = \underset{x \in X} \to \bigcup{\A(T_x X)}\arrow X.
$$
We denote by $\E(=\E(X))$ the set of $C^\infty$ global sections of $\A(X)$, 
$$
\E(X) = \Gam (X, \A(X)).
$$
Let $\Phi^0$ be a closed element of $\E$. 
Then we have the vector spaces $E^k
(T_xX)$  for each $x \in X$ and $k \geq 0$.
 We define the vector bundle $E^k_{\ss-style X}(=E^k)$ over $X$  as
$$
E^k_{\ss-style X}: = \underset {x \in X}\to\bigcup{E^k (T_x X)}\arrow X.
$$ for each $k \geq 0$.
(Note that the  fibre of $E^1$ is $\frak g/\frak h$.)
Then we define the graded module $\Gam (E)$ over $\Gam (\w^*)$ as $\oplus_k
\Gam (E^k)$, where $\Gam$ denotes the set of global C$^\infty$ sections and 
$\w^p$ is the sheaf of germs of smooth $p$ forms on $X$.
\proclaim {Theorem 1-3 } 
$\Gam (E)$ is the differential graded module in 
$\oplus_k\Gam ( \oplus_i \w^{ p_i +k -1} )$ with respect to the exterior
derivative $d$. 
\endproclaim
\demo{Proof}
Since $\Gam (E)$ is the graded module generated by $\Gam (E^0)$,
it is suffices to prove that $di_v \Phi^0 $ is an element of $\Gam (E^1)$ for 
$v \in \Gam (TX)$.
We denote by Diff$(X)$ the group of diffeomorphisms of $X$. 
Then there is the action of Diff$(X)$ on differential forms on $X$ and 
we see that  $\E(X)$ is invariant under the action of Diff$(X)$. 
An element of $\Gam(E^0)$ is given as $i_v \Phi^0= (i_v \phi_1, \cdots, i_v
\phi_l)$,  where $v \in \Gam (TX)$. 
Since $\Phi^0$ is closed, we have 
$$
d i_v \Phi^0 = L_v \Phi^0.
$$
The vector field $v$ generates the one parameter group of transformation $f_t$. 
Then $L_v \Phi^0 = \frac d{dt}f^*_t \Phi^0 |_{t=0}$. 
Since $\E(X)$ is invariant under the ation of Diff$(X)$, 
$f^*_t ( \Phi^0) \in \E(X)$.
Since the tangent space of $\E$ at $\Phi^0$ is $\Gam(E^1)$, 
$L_v \Phi^0 \in \Gam(E^1)$. 
Hence $d i_v\Phi^0 \in \Gam(E^1)$. 
From definition of $E^k(V)$, we see that 
$da \in \Gam (E^k)$ for all $a \in \Gam(E^{k-1})$ for all $k$.
\qed\enddemo
Then from theorem 1-3, 
we have a complex $\#_{\Phi^0}$ 
$$
\CD
\Gam (E^0) @>d_0>>\Gam (E^1)@>d_1>>\Gam(E^2)@>d_2>>\cdots,
\endCD
 \tag{\#$_{\Phi^0}$}
$$
 where $\Gam (E^i)$ is the set of $C^\infty$ global sections for each vector
bundle and 
 $d_i = d|_{E^i}$ for each $i =0,1,2$.
The complex \#$_{\Phi^0}$ is a subcomplex of 
the direct sum of the de Rham complex (For simplicity, we call this
complex the de Rham complex):
$$
\CD 
@.\Gam (E^0)@>d_0>>\Gam (E^1)@>d_1>>\Gam (E^2)@>d_2>>\cdots \\ 
@.@VVV                    @VVV         @VVV \\
\cdots@>d>>\Gam (\oplus_i \w^{p_i-1} ) @>d>> \Gam ( \oplus_i \w ^{p_i} ) @>d>>\Gam
(\oplus_i
\w^{p_i+1} )@>d>>\cdots. 
\endCD
$$
If $\O$ is an elliptic orbit, the complex $\#_{\Phi^0}$ is an elliptic complex 
for all closed $\Phi^0 \in \E^1$ on any $n$ dimensional compact manifold $X$
( Note that the complex in definition 1-1 is the symbol complex of
$\#_{\Phi^0}$ ). Then we have a finite dimensional cohomology group
$H^k(\#_{\Phi^0})$ of  the elliptic complex $\#_{\Phi^0}$.
Since $\#_{\Phi^0}$ is a subcomplex of deRham complex, there is the map $p^k$ from
the cohomology group of the complex
$\#_{\Phi^0}$ to de Rham cohomology group: 
$$
p^k\:H^k (\#_{\Phi^0} )  \arrow \underset i\to\oplus H^{p_i-k+1}(X,\R).
$$ 
where
$$
H^k ( \#_{\Phi^0} ) = \{ \, \a \in \Gam (E^k) \, |\, d_k \a =0 \, \} /
\{\, d\b \,| \, \b \in \Gam ( E^{k-1})\, \}.
$$
\proclaim{Definition 1-4 (Topological calibrations and topological orbits) } 
A closed element $\Phi^0 \in \E(X)$ is a topological calibration if the map 
$$
p^k\:H^k (\#_{\Phi^0} )  \arrow \underset i\to\oplus H^{p_i+k-1}(X,\R)
$$
is injective for $k=1,2$.
An orbit $\O$ in $\oplus_i \w^{p_i }V^*$  is topological over a manifold $X$
if any closed element of $\E(X) $ is a topological calibration. 
An orbit $\Cal O$ is topological if $\Cal O$ is topological
over any compact 
$n$ dimensional manifold $X$.
\endproclaim
\proclaim{Lemma 1-5}
Let $\Cal O$ be a metrical orbit and $\Phi^0$ an element of 
$\E = \Gam (X, \A_{_\O}(X) )$.
Then there is a canonical metric $g_{\ss-style\Phi^0}$on $X$ corresponding to 
each $\Phi^0 $. 
\endproclaim
\demo{Proof}
The orbit $\O$ is defined in terms of $\Phi^0_V \in \oplus_i
\w^{p_i}V^*$  on $V$. We also have $\Phi^0 (x)\in \A(T_x X)$ on each tangent
space $T_x X$.  Let $\text{Isom} (V, T_xX ) $ be the set of isomorphisms
between 
$V$ and $T_xX$.
Then define $H_x$ by 
$$
H_x = \{  h \in \text{Isom} (V, T_xX ) \, |\,  \Phi^0_V = h^*\Phi^0(x) \, \}.
$$
Then we see that $H_x$ is isomorphic to the isotropy group $H$. 
$h_* g_V$ defines the metric on the tangent space $T_x X$ for $h \in
H_x$. Since $H$ is a subgroup of O$(V)$, the metric $h_*g_V$ does not depend
on a choice of $h \in H_x$.
\qed\enddemo
Let $\O$ be an orbit in  $\oplus_i \w^{p_i }V^*$ . Then 
we define the moduli space $\M_{_\O} (X)$ by 
$$
\M_{_\O} (X) = \{ \, \Phi \in \E \, |\, d\Phi =0 \, \} / \text{Diff}_0(X),
$$
where Diff$_0(X)$ is the identity component of the group of diffeomorphisms
for $X$.  We denote by $\widetilde{\M}_{_\O} (X)$ the set of closed elements in
$\E$: 
$$
\wtil{\M}_{\O}(X) =\{\, \Phi \in \E(X) \, |\, d\Phi =0 \,\}.
$$  We have the natural projection $\pi \: \widetilde{\M}_{_\O} (X) \to \M_{_\O}
(X)$.
Let $\Phi^0$ be an element of $\widetilde{\M}_{_\O} (X)$. 
As we shall show in section 2, $\E(X)$ is regarded as
 a infinite dimensional homogenous space (a Hilbert manifold). Hence we have the
tangent space 
$T_{\ss-style\Phi^0}\E(X)$. We denote by $\Cal H$ the Hilbert space consisting of 
closed forms. Then the space $\wtil{\M}_{\O}(X)$ is the intersection 
between the Hilbert space $\Cal H$ and the Hilbert manifold $\E(X)$. 
We define an infinitesimal tangent space  of $\wtil{\M}_\O$ by the intersection 
$\Cal H\cap T_{\ss-style\Phi^0}\E$.
Since $T_{\ss-style\Phi^0}\E(X)=E^1$, the infinitesimal tangent space is written as 
$$
\Cal H\cap T_{\ss-style\Phi^0}\E(X) = \Cal H \cap E^1. 
$$  Then we shall discuss if the infinitesimal tangent space
is regarded as the tangent space of actual deformations.
\proclaim{Definition 1-6 } 
A closed element $\Phi^0 \in \E(X)$ is unobstructed if there 
exists an integral curve $\Phi_t(a)$ in $\wtil{\M}(X)$ for each $a \in\Cal H \cap
E^1$ such  that 
$$
\frac d{dt}\Phi_t(a)|_{t=0} = a
$$
An orbit $\O$ is unobstructed if any $\Phi^0 \in \wtil{\M}_{\O}(X)$  is unobstructed 
for any compact $n$ dimensional manifold $X$.
(see  section 2 for the precise statement with respect to Sobolev norms.)
\endproclaim
The following figures (i),(ii) and (iii) explain our situation well.
If the Hilbert space $\Cal H$ is in a generic position, 
the intersection $\wtil{\M}(X)=
\E\cap \Cal H$ is  smooth and every element $\Phi^0$ is unobstructed (see figure (i)). 
However if $\Cal H$ is in a special position, $\wtil{\M}(X)$ may be singular. 
(see figure (ii)). Then the infinitesimal tangent space may not gives 
actual deformations. In figure (ii), the point $p$ is singular and 
the infinitesimal tangent space at $p$ coincides with $\Cal H$. 
Figure (iii) shows this problem is subtle. The intersection $\E$ is a line and 
it seems to be non-singular. However the infinitesimal tangent space is $\Cal H$ at each
point, which is obstructed.
\bigskip
\bigskip\bigskip
\bigskip\bigskip
\bigskip\bigskip
\bigskip
\bigskip
\bigskip
\bigskip
\bigskip
\bigskip
figure(i)\qquad\qquad\qquad\qquad figure(ii)\qquad\qquad\qquad\qquad figure(iii)
\bigskip
\bigskip
We shall prove the following theorems in section 2.
\proclaim{Theorem 1-7 ( Criterion of unobstructedness)} 
We assume that an orbit $\O$ is elliptic. 
If the map $p^2\:H^2(\#_{\ss-style\Phi^0})\to \oplus_i H^{p_i+1}_{\ss-style DR}(X)$
is injective, then $\Phi^0$ is unobstructed.
\endproclaim
We shall prove the following theorems in section 3.
\proclaim{Theorem 1-8 } 
If an orbit $\O$ is metrical, elliptic and topological, then the corresponding moduli space $\M_{_\O} (X)$ is a smooth manifold. 
( In particular $\M_{_\O} (X)$ is Hausdorff .)
Further $\M_{_\O} (X)$ has canonical coordinates given by an open ball of the
cohomology group H$^1 (\#_{\Phi}) $. 
\endproclaim
Since de Rham cohomology group is invariant under the action of Diff$_0(X)$, we
have the map 
$$
P\: \M_{_\O} (X) \arrow \underset i \to \oplus H^{p_i}_{dR} (X).
$$
Then we have 
\proclaim{Theorem 1-9}
If an orbit $\O$ is metrical, elliptic and topological, then the map $P$ is locally injective. 
\endproclaim
Further under the assumption that $\O $ is metrical, elliptic and topological 
we have
\proclaim{Theorem 1-10} 
Let $I(\Phi$) be the isotropy group, 
$$
I(\Phi) = \{ \, f \in \text{\rm Diff}_0(X) \, |\,  f^* \Phi = \Phi \, \}.
$$
Then there is a sufficiently small  slice $S_{\Phi^0}$ at $\Phi^0 $
such that the isotropy group $I(\Phi^0)$ is a subgroup of $I(\Phi)$ for each
$\Phi\in S_{\Phi^0}$,i.e.,
$$
I(\Phi^0) \subset I(\Phi).
$$
(Our definition of the slice will be given in section 2 and 3.)
\endproclaim
\proclaim{Theorem 1-11}
Let $\wtil{\M}_{\ss-style\O}(X)$ be the set of closed elements of $\E^1$. 
We denote by $\Diff(X)$ the group of diffeomorphism of $X$. There is the
action of $\Diff(X)$ on $ \wtil{\M}_{\ss-style\O}(X)$.  Then the quotient
$\wtil{\M}_{\ss-style\O}(X)/\Diff(X)$ is an orbifold.  
\endproclaim

\head 
\S2.  Deformations of calibrations
\endhead
\subhead \S 2-0 Preliminary results
\endsubhead
Let $X$ be a manifold and we denote by  $\w^*$ the 
differential forms on $X$. 
Let $P$ be a linear operator acting on $\w^*$. 
Then the operator $P\: \w^* \to \w^*$ is a derivative if 
$P$ satisfies the followings:
$$\align 
&P( s+t ) =P(s) + P(t),\\
&P(s\w t) = P(s) \w t +s\w P(t),
\endalign
$$
where $s,t \in \w^*$.
An anti-derivative  $Q$ is also a linear operator defined by 
the following: 
$$
\align 
&Q( s+t) = Q(s) + Q(t), \\
&Q(s\w t) = Q(s)\w t + (-1)^{|s|}s\w Q (t),
\endalign$$
where $|s|$ denotes the degree  of a differential form $s$.
Then the exterior derivative $d$ is the anti-derivative and  
the  differential representation $\h{\rho}_a$  is a derivative 
for each $a\in $End$(TX)$.
\proclaim{Lemma 2-0-1} 
The commutator $[\h{\rho}_a, d]=\h{\rho}_a\circ d -d\circ\h{\rho}_a$ is the anti-derivative. 
We denote by $L_a$ the commutator $[\h{\rho}_a, d]$.
\endproclaim
\demo{Proof} 
In general the commutator of a derivative $P$ and 
an anti-derivative $Q$ is an anti-derivative if $Q$ preserves degrees of differential forms.
\enddemo
The operator $L_a$ is regarded as a generalizations of the Lie 
derivative. Indeed we have
\proclaim{Lemma 2-0-2} The commutator $L_a$ is expressed as 
$$
L_a \:\w^n \arrow \w^{n+1},
$$
$$\align
L_a \eta (u_0,u_1, \cdots, u_n)= &\sum_{i=0}^n (-1)^i L_{a\,u_i}\eta\, 
( u_0, u_1,\overset \check{i}\to{\cdots}, u_n),\\
-&\sum_{i<j}(-1)^{i+j}\eta( a[u_i,u_j ],u_0,\overset{\check{i}\quad\check{j}}\to{\cdots\cdots},u_n ) 
\endalign$$
where $\eta$ is an $n$ form and $a\in $End$(TX)$ maps a vector $u_i $ to 
$au_i \in TX$ and we denote by $L_{au_i}$ the ordinary Lie derivative. 
\endproclaim
\demo{Proof} 
It is sufficient to show the lemma with respect to vectors $\{u_i\}$ satisfying 
$[u_i, u_j ]=0$.
Then we have 
$$
\align 
&(\hrho_a d\eta )(u_0, \cdots u_n )  = \sum_i (-1)^i (i_{au_i} d\eta) ( u_0,
\overset\check{i}\to{
\cdots},u_n ) \\
&(d\hrho_a \eta)(u_0, \cdots , u_n ) = - \sum_i (-1)^i (d i_{au_i} \eta )(u_0,
\overset\check{i}\to{
\cdots},u_n ) .
\endalign
$$ 
Hence from $L_{au_i} = d i_{au_i} + i_{au_i} d$, we have the result.
\enddemo
we also have a description of the commutator between 
$L_a$ and $\h{\rho}_a$,
\proclaim{Lemma 2-0-3} 
$$
[L_a, \h{\rho}_b] = i_{N(a,b)} - L_{ab},
$$
where $a,b \in$\text{\rm End}$(TX)\cong \w^1\otimes  T$ and 
a tensor $N(a,b)\in  \w^2\otimes T$ is given by the following
$$\align
N(a,b) (u,v) = &ab [u,v]+ba[u,v]+[au,bv]-[av,bu]\\ -&a[bu,v] +a[bv,u]-b[au,v]+b[av,u],
\endalign$$
for $u,v\in TX$, 
and  $i_{N(a,a)}$ is the composition of the interior product and the wedge product of the tensor
$N(a,a)\in \w^2\otimes TX$. 
\endproclaim
\demo{Remark} 
The tensor $N(a,b)$ is a generalization of the Nijenhuis tensor.
\enddemo
\demo{Proof of lemma 2-0-3} 
For $a,b \in $End$(TX)$, we have the tensor $N(a,b)\in \w^2\otimes TX$.
Then $i_{N(a,b)}$ is the linear operator from $\w^* \to \w^{*+1}$. 
We see that $i_{N(a,b)}$ is an anti-derivative. 
By lemma 2-0-1, $L_{ab}$ is an anti-derivative, where $ab$ denotes the composition of
endmorphisms. As in proof of lemma 2-0-1, the commutator $[L_a , \hrho_a ]$ is also an
anti-derivative.  
Hence it sufficient to show that the identity in lemma 2-0-3 for functions and $1$
forms.  For a function $f$, we have 
$[L_a,\h{\rho}_a]f = -\h{\rho}_a L_a f = - L_{a^2}f$. 
Since $i_{N(a,b)} f =0$, we have the identity. 
For a one form $\theta$ by lemma 2-0-2, we have 
$$\align
L_a \h{\rho}_b \theta ( u,v) = &(L_{au}\h{\rho}_b\theta)(v)-( L_{av}\h{\rho}_b\theta)(u) +\h{\rho}_b\theta(\h{\rho}_a[u,v]) \\
=&au(\h{\rho}_b\theta (v) )-av(\h{\rho}_a\theta(u)) +\theta( ba[u,v])\\
-&\h{\rho}_b\theta([au,v]) + \h{\rho}_b\theta([av,u]).\\
\\
\\
\h{\rho}_bL_a\theta (u,v) =& 
(L_a\theta)(\h{\rho}_b u, v) + (L_a\theta)(u,\h{\rho}_bv) \\
=&(L_{abu}\theta)(v)-(L_{av}\theta)(\h{\rho}_b u) +\theta(a[bu,v]) \\
+&(L_{au}\theta)(\h{\rho}_bv)-(L_{abv}\theta)(u)+\theta(a[u,bv])\\
=&(abu)(\theta v)-(av)(\theta (bu))+\theta(a[bu,v]) \\
-&\theta([abu,v])+\theta([av,bu])+\theta(a[u,bv])\\
+&(au)\theta(bv)-(abv)\theta(u)\\
-&\theta([au,bv])+\theta([abv,u])
\endalign$$
Hence the commutator is given by  
$$
\align 
[L_a,\h{\rho}_b]\theta(u,v)=
&-(abu)(\theta v)+\theta([abu,v])+(abv)\theta(u) -\theta([abv,u])\\ 
+&\theta( ba[u,v])+\theta([au,bv])-\theta([av,bu])\\
-&\theta(a[bu,v])+\theta(a[bv,u])-\theta(b[au,v]) + \theta(b[av,u])\\
=&-L_{ab}\theta(u,v)+i_{N(a,b)}\theta
\endalign
$$
\qed\enddemo
\proclaim{Lemma 2-0-4}
We assume that  $\Phi$ and $\h{\rho}_a\Phi$ are closed forms 
respectively.  Then
$d\h{\rho}_a\h{\rho}_a \Phi$ is an element of $\Gam (E^2)$.
\endproclaim
\demo{Proof} 
$$\align
d\h{\rho}_a \h{\rho}_a \Phi =& \h{\rho}_a d\h{\rho}_a\Phi -L_a\h{\rho}_a\Phi = -L_a\h{\rho}_a\Phi\\
=&-\h{\rho}_a L_a \Phi -i_{N(a,a)}\Phi + L_{a^2}\phi.
\endalign$$
Since $L_a\Phi =\h{\rho}_ad\Phi -d\h{\rho}_a\Phi =0,$
we have 
$$
d\h{\rho}_a \h{\rho}_a \Phi= -i_{N(a,a)}\Phi + L_{a^2}\Phi.
$$
Since $N(a,a) \in \w^2\otimes T \cong \w^1\otimes $End$(TX)$, 
then it follows from our definition of $E^2$ that 
$$
i_{N(a,a)}\Phi \in \Gam(E^2).
$$
Since $L_{a^2}\Phi = -d\h{\rho}_{a^2}\Phi\in d\Gam(E^1)\subset \Gam(E^2)$.
Hence we have the result.
\qed\enddemo 
We denote by $G=G(a,a)$ the operator 
$i_{N(a,a)}-L_{a^2}$. Then we consider the  commutator 
$[\hrho_a, G(a,a) ]$.  For simplicity we write this by 
$Ad_{\hrho_a}G(a,a) (= Ad_{\hrho_a}G )$,
$$
Ad_{\hrho_a}G(a,a) = [\hrho_a, G(a,a)].
$$ 
The $k$th composition of commutator 
is denoted by 
$$
Ad_{\hrho_a}^k G = [\hrho_a ,[\hrho_a, \cdots [ \hrho_a , G(a,a) ],\cdots]],
$$
where Ad$_{\hrho_a}G(a,a)$ acts on differential forms.
\proclaim{Lemma 2-0-5} 
$Ad_{\hrho_a}^kG(a,a)\Phi^0$ is an element of $\Gam(E^2)$.
\endproclaim 
\demo{Proof} 
At first we consider Ad$_{\hrho_a}G(a,a)\Phi^0$. 
By lemma 2-0-3, we have 
$$\align
Ad_{\hrho_a}G(a,a)\Phi^0 =& [ \hrho_a, G(a,a) ]\Phi^0 \\
=& [\hrho_a , i_{\ss-style {N(a,a)}}]\Phi^0 - 
[\hrho_a , L_{a^2} ]\Phi^0 \\
=& [\hrho_a , i_{\ss-style {N(a,a)}}]\Phi^0 
+G(a^2 ,a ) \Phi^0.
\endalign
$$
Since $N(a^2,a)\in \w^1\otimes $End$(TX)$,
as in lemma 2-0-4, $G(a^2,a) \Phi^0$ is an element of $\Gam(E^2)$. 
We see that 
$[\hrho_a, i_{\ss-style{N(a,a)}}]$ is give by the 
interior product of the tensor $\hrho_a (N(a,a)) \in \w^1\otimes $End$(TX)$, where 
$\hrho_a$ acts on the tensor $N(a,a)$.
Hence $[\hrho_a, i_{\ss-style{N(a,a)}}]\Phi^0$ 
is an element of  $\Gam(E^2)$. 
Therefore $Ad_{\hrho_a}G(a,a)\Phi^0 \in \Gam(E^2)$. 
By induction, we see that 
Ad$_{\hrho_a}^kG(a,a) \Phi^0$ is an element of  $\Gam(E^2)$.
\qed\enddemo

\subhead \S2-1 Primary obstruction
\endsubhead 
In this section we use the same notation as in section 1 and subsection
2-0.  The references for analytic tools are found in [7],[9],[18],and [19]. 
Let $X$ be a real $n$ dimensional compact manifold. We fix a Riemannian
metric $g$ on $X$. 
(Note that this metric does not depend on calibration $\Phi$.) We denote by
$C^\infty(X,\w^p)$ the set of smooth
$p$ forms on
$X$. Let $L^2_s ( X, \w^p)$ be the Sobolev space and suppose that $s > k +
\frac n2$., i.e., the completion  of $C^\infty(X, \w^p)$ with respect to the
Sobolev norm $\|\, \|_s$, where $k$ is sufficiently large ( see [9] for
instance).  Then we have the inclusion $ L^2_s ( X, \w^p ) \arrow C^k ( X,
\w^n )$.  We define $\E_s$ by 
$$
\E_s = C^k ( X, \A_\O(X) ) \cap L^2_s ( X, \oplus_{i=1}^l\w^{p_i} ).
\tag2-1-1$$
Then we have 
\proclaim{Lemma 2-1-1}
$\E_s$ is a Hilbert manifold (see [19] for Hilbert manifolds ). The
tangent space $T_{\Phi ^0}\E_s$ at
$\Phi^0$ is given by 
$$
T_{\Phi^0} \E_s = L^2_s ( X, E^1).
$$ 
\endproclaim
\demo{Proof}
We denote by exp the exponential map of Lie group $G=$GL$(n,\R)$. Then 
we have the map $k_x$ 
$$
k_x \: E^1 ( T_xX ) \arrow \A ( T_x X),
\tag2-1-2$$
by 
$$
k_x ( \hrho_\xi\Phi^0(x) ) = \rho_{\exp \xi}\Phi^0 (x).
\tag2-1-3$$
for each tangent space $T_x X$.
From 2-1-2, we have the map $k$
$$
k \:L^2_s ( E^1) \arrow \E_s,
\tag2-1-4$$by 
$$
k|_{E^1(T_x X)} = k_x.
$$
The map  $k$ defines local coordinates of $\E_s$.
\qed\enddemo
Let 
GL$(TX)$ be the group of gauge transformations, i.e., 
for $g \in $GL$(TX)$ we have the diagram:
$$
\CD
TX @>g>> TX \\
@VVV@VVV \\
X@>id>>X
\endCD
$$
An element $g\in $GL$(TX)$ acts on $\E_\O(X)$
by 
$$
\Phi \mapsto \rho_g ( \Phi)
$$
The tangent space $T_{\Phi^0}\E(X)$ is given by $E^1(X)$, 
$$
E^1(X) =\{\, \h{\rho}_a \Phi^0 \, |\, a \in End(TX) \, \}
$$
where $\h{\rho}$ is the differential representation of $\rho$.
We denote by H$(TX)$ be the gauge transformations with structure group 
H, i.e., the isotropy group. 
Then by lemma 2-1-1, $\E$ is regarded as the infinite dimensional homogenous
space  GL$(TX)/$H$(TX)$. 
Let $\Cal H$ be the closed subspace of $L^2_s ( X, \oplus_{i=1}^l\w^{p_i} )$ 
consisting of closed forms. 
Then $\wtil{\M}_s(X)$ is the intersection between $\E$ and $\Cal H$. 
The image $dE^0(X)$ is given by 
$$
dE^0(X)=\{ \, d i_v \Phi^0 = L_v \Phi^0\, |\, v \in TX\,\},
$$
where $L_v $ is the Lie derivative with respect to $v \in TX$. 
Hence the cohomology H$^1( \#)$  of the complex $\#_{\Phi^0}$ is considered as the infinitesimal tangent space of 
the moduli space $\M(X)=\wtil{\M}(X)/\Diff(X)$.
However, the moduli space may not be a manifold in general. 
This is because the infinitesimal tangent space may not be exponentiate the actual deformations. Then there exists an obstruction. 
This is a general problem of deformation. 
In our situation, we must study the intersection $\E\cap \Cal H$.  In order
to obtain a deformation space,  we shall construct a deformation of $\Phi^0$
in terms of a power series in $t$. We consider a formal power
series in $t$: 
$$
a(t)= a_1 t+ \frac1{2!}a_2 t^2 +\frac1{3!}a_3 t^3 + \cdots \in End(TX)
[[t]],
\tag 2-1-5$$
where $a_k \in $ End $(TX)$.
We define a formal power series $g(t)$ by,
$$
g(t) = \exp a(t) \in GL(TX)[[t]]
$$
For simplicity, we put $a=a(t).$
The gauge group GL$(TX)$ acts on differential forms by $\rho$.
This action
$\rho$ is written in terms of the differential representation $\hrho$, 
$$\align
\rho_{g(t)} \Phi^0 = &\Phi^0 +\hrho_a\Phi^0 + 
\frac1{2!}\hrho_a\hrho_a \Phi^0 +
\frac1{3!}\hrho_a\hrho_a\hrho_a\Phi^0+\cdots\\ =&\Phi^0 + \h{\rho}_{a_1}
\Phi^0 t +  
\frac12 (\h{\rho}_{a_2}\Phi^0+\h{\rho}_{a_1}\h{\rho}_{a_1}\Phi^0  )t^2 + \cdots ,
\tag 2-1-6\endalign
$$
where $\hrho$ is just written as
$$
\hrho_{a(t)}\Phi^0 = \sum_{k=1}^\infty \frac1{k!}\hrho_{a_k}\Phi^0 t^k.
$$
The equation what we want to solve is , 
$$
d \rho_{g(t)}\Phi^0 =0.
\tag eq$_*$
$$
We must find a power series $a=a(t)$ satisfying the condition (eq$_*$). 
At first we take $a_1$ such that $d\hrho_{a_1}\Phi^0 =0$. 
Then it remains to determine $a_2, a_3,\cdots$ satisfying (eq$_*$). 
$d\rho_{g(t)}\Phi^0$ is written as a power series, 
$$
\align
d\rho_{g(t)}\Phi^0 =&\sum_{k=1} \frac1{k!}dR_k t^k,\tag 2-1-7\\
\endalign
$$
where $R_k$ denotes the homogenous part of degree $k$.
Hence the equality $d\rho_{g(t)}\Phi^0=0$ is reduced to the system of 
infinitely many equations
$$
dR_k=0, \quad  k=1,2,\cdots
\tag eq$_{k}$
$$
By our assumption  
$
d \h{\rho}_{a_1}\Phi^0 =0, $ we already have $dR_1 =0$ (see (2-1-6)). 
Thus in order to obtain $a(t)$, it suffices to determine $a_{k}$ 
satisfying (eq$_{k}$) by induction on $k$. 
By (2-1-6), the term of the second order $dR_2$ is given as
$$
dR_2 =\frac1{2!} \( d\h{\rho}_{a_2} \Phi^0 +  
d\h{\rho}_{a_1}\h{\rho}_{a_1}\Phi^0  \) 
\tag2-1-8$$
We denote by $Ob_2(a_1)$ the quadratic term, 
$$
Ob_2 (a_1)= \frac1{2!} (d\h{\rho}_{a_1}\h{\rho}_{a_1}\Phi^0 )
\tag2-1-9$$
Then by lemma 2-0-4 in section 2-0, $Ob_2$ is an element of $\Gam(E^2)$,
which is explicitly written as
$$
Ob_2(a_1) =-\frac1{2!} (-i_{N(a_1, a_1)} +L_{a_1^2}  )\Phi^0,
\tag2-1-10
$$
Since $Ob_2(a_1)$ is a d-closed form, this defines a representative of the
cohomology group H$^2(\#)$.  In order to determine $a_2$ satisfying
$dR_2=0$, we must solve the equation,
$$
\frac1{2!}d \hrho_{a_2} \Phi^0= -Ob_2(a_1).
\tag {$eq_2$}
$$
The L.H.S of $(eq_2)$ cohomologically vanishes in H$^2(\#)$.
 Hence if the class $[Ob_2(a_1)]\in H^2(\#)$ does not vanishes, 
 there exists no solution $a_2$ of $eq_2$ and no deformation with 
 $a_1$. In this sense we call the class $[Ob_2(a_1)]$ 
the obstruction to deformation of $\Phi^0$ ( the primary obstruction ).
 If $[Ob_2(a_1)]$ vanishes, then we have a solution $a_2$ by 
 $$
 \frac1{2!}\hrho_{a_2}\Phi^0 = -d_1^* G_\# (Ob_2 (a_1) ),
 \tag2-1-11$$
where $G_\#$ denotes the Green operator of the complex $\#$.
 It is quite remarkable that the representative $Ob_2(a_1)$ is d-exact form. 
 Hence $Ob_2(a)$ is in kernel of the map $p^2\: H^2(\#) \to 
 \oplus_i H^{p_i+1}(X)$.  Hence we obtain a nice criterion of unobstructedness.
 \proclaim{Theorem 2-1-2} 
If the map $p^2 \: H^2(\#) \to \oplus_i H^{p_i +1}(X)$ is injective, 
The obstruction class $[Ob_2(a_1)] $ vanishes.
\endproclaim
\head \S 2-2 Higher obstructions 
\endhead
Similarly we obtain infinitely many obstructions to deformation of $\Phi^0$. We define an operator $G(a,a)$ on $\w^*$ by 
 $$
 G(a,a) = i_{N(a,a)} -L_{a^2},
 \tag2-2-1$$
 where $a=a(t)\in$End$(TX)[t]$. We denote its $k$th homogenous part 
 by $G(a,a)_k$. Then by lemma 2-0-4, we have 
 $$
 Ob_2 (a_1) = -\frac1{2!}G(a,a)_2.
 \tag2-2-2$$
 We assume that $a_1, a_2,\cdots a_{k-1}$ are determined satisfying 
 $dR_1=0, dR_2=0,\cdots ,dR_{k-1}=0$. 
 Then $dR_k$ is written as a $d$-exact form:
$$
dR_k = d\hrho_{a_k} \Phi^0 + \sum_{l=2}^k \frac1{l!}(d\hrho_a^l )_k\Phi^0,
\tag 2-2-3$$
where $(d\hrho_a^l)_k$ denotes the $k$ th homogeneous part of
$d\hrho_a^l$. We define $Ob_k(a_{\ss-style{<k}} )$ as $\sum_{l=2}^k
\frac1{l!}(d\hrho_a^l)_k
\Phi^0,$ 
where $a_{<k}= a_1 t + \frac1{2!} a_2t^2 + \cdots +\frac1{(k-1)!}a_{k-1}
t^{k-1}.$ 
Then we have
\proclaim{Proposition 2-2-1}
 $$
 dR_k =\frac1{k!}d\hrho_{a_k}\Phi^0 +Ob_k( a_{<k} ),
 $$where $Ob_k$ is written as
 $$\align
&{\s-style{ 
 Ob_k(a) = \Big(-\frac1{2!}G(a,a) \Phi^0 +\frac1{3!} [\hrho_a, G(a,a) ]
 \Phi^0- \cdots 
 + (-1)^{k-1} \frac1{k!}[\hrho_a, [\cdots,[\hrho_a, G(a,a) ]]\cdots ]
\Big)_k \Phi^0 }}\\
 =&\(f(Ad_{\hrho_a})G(a,a)\)_k \Phi^0,
 \endalign
 $$
 where $f(x) $ is a convergent sequence, 
 $$
 f(x)= -\frac1{2!} + \frac1{3!}x-\frac1{4!} x^2 -\cdots=-\frac{e^{-x} -1+x}{x^2}
 \tag2-2-3$$
 and Ad$_{\hrho_a}$ is the adjoint operator $[\hrho_a, \, ]$.
 Substituting $Ad_{\hrho_a}$ into $f(x)$, we have an operator 
 $f(Ad_{\hrho_a})$.  The higher obstruction $Ob_k$ consists of commutators.
Hence 
 $Ob_k =(f(Ad_{\hrho_a})G(a,a))_k\Phi^0$ is essentially the interior
product of 
 $\Phi^0$ in terms of the tensors of type $ \w^2\otimes T$.
 Hence we see that $Ob_k(a_{<k}) \in E^2$.
 \endproclaim
\demo{Proof} 
In the case $k=1$ we have the proposition. 
We shall prove the proposition by induction on $k$. 
We assume that proposition holds for all $l < k$. 
Then we have 
$$
dR_l =-(L_a )_{l}\Phi^0+\( f(\text{Ad}_{\hrho_a})G(a,a)\)_{l}\Phi^0.
\tag2-2-4$$
We put $(L_a)_{\ss-style{<k}}$ as 
$$
(L_a)_{\ss-style{<k}}= \sum_{l=2}^{k-1} (L_a)_l.
$$
If $dR_l=0\, ( l < k)$, from our assumption 
we have 
$$\align
(L_a)_{\ss-style{<k}} \Phi^0 =& -\( f(\text{Ad}_{\hrho_a})G(a,a)\)_{\ss-style{<k}}\Phi^0\\
=&\(-\frac1{2!}G(a,a)+\frac1{3!}[\hrho_a, G(a,a)]-\frac1{4!}[\hrho_a,[\hrho_a,G(a,a)]]+\cdots \)_{\ss-style{<k}}\Phi^0.\\
=&\sum_{l=2}^k (-1)^{l-1}\frac1{l!}
(\text{Ad}_{\hrho_a}^{l-2}G(a,a))_{\ss-style{<k}} \Phi^0
\tag 2-2-5\endalign
$$
Then by using lemma 2-0-3, we have 
$$\align
d(\rho_{e^a})_k&\Phi^0 =\sum_{l=1}^k\frac1{l!}(d\hrho_a^l )_k\Phi^0\\
=&-(L_a)_k\Phi^0-\frac1{2!}(G(a,a)+ 2\hrho_a L_a)_k \Phi^0 \\
&- \frac1{3!}(G(a,a)\hrho_a+2 \hrho_aG(a,a) +3\hrho_a\hrho_aL_a)_k\Phi^0\\
&-\frac1{4!}(G(a,a)\hrho_a\hrho_a+2\hrho_aG(a,a)\hrho_a+3\hrho_a\hrho_aG(a,a)+ 
4\hrho_a\hrho_a\hrho_aL_a )_k\Phi^0-\cdots. \\
\tag 2-2-6\endalign 
$$
Since the degree of $a=a(t)$ is greater than or equal to one,  we have 
$$
( \hrho_a^m L_a )_k =  ( \hrho_a^m (L_a)_{\ss-style{<k}} )_k,
\tag2-2-7$$
for a positive integer $m$.
Hence from (2-2-7), we substitute (2-2-5) into (2-2-6) and we have
$$\align
&d(\rho_{e^a})_k\Phi^0=-(L_a)_k\Phi^0-\frac1{2!}G(a,a)_k \Phi^0\\
&-\frac1{2!}2(\hrho_a(-\frac1{2!}G(a,a)+\frac1{3!}\text{Ad}_{\hrho_a}G(a,a)+\cdots))_k\Phi^0\\
&-\frac1{3!}(G(a,a)\hrho_a+2\hrho_a G(a,a))_k \Phi^0-\frac1{3!} 3(\hrho_a\hrho_a (-\frac1{2!}G(a,a)+\cdots))_k\Phi^0-\cdots. \\
\endalign $$
Then we calculate each homogeneous part with respect to $a$ and 
we have 
$$\align
&d(\rho_{e^a})_k\Phi^0=-(L_a)_k\Phi^0-\frac1{2!}G(a,a)_k\Phi^0\\
&+
(\frac{2}{2!2!}\hrho_a G(a,a)-\frac2{3!}\hrho_a G(a,a)-\frac1{3!}G(a,a)\hrho_a )_k\Phi^0\\
&+(-\frac{2}{2!3!}\hrho_a[\hrho_a, G(a,a)]   +
\frac3{3!2!}\hrho_a\hrho_aG(a,a) )_k\Phi^0 \\
&+\frac1{4!}(-G(a,a)\hrho_a\hrho_a -2\hrho_a GG(a,a) \hrho_a -3\hrho_a\hrho_aG(a,a)))_k\Phi^0 + \cdots \\
=&-(L_a)_k\Phi^0-\frac1{2!}G(a,a)_k\Phi^0+\frac1{3!}[\hrho_a, G(a,a)]_k\Phi^0\\
&+(-\frac1{4!}G(a,a)\hrho_a\hrho_a+
(\frac{-2}{4!}+\frac2{3!2!})\hrho_aG\hrho_a+(\frac{-3}{4!}+\frac3{3!2!}-\frac2{3!2!})\hrho_a\hrho_a G )_k \Phi^0+\cdots\\
=&-(L_a)_k\Phi^0-\frac1{2!}G(a,a)_k\Phi^0+\frac1{3!}[\hrho_a, G(a,a)]_k\Phi^0-\frac1{4!}[\hrho_a,[\hrho_a, G(a,a)]]_k\Phi^0 +\cdots\\
=&-(L_a)_k\Phi^0 +\sum_{l=2}^k(-1)^{l-1}\frac1{l!}\text{Ad}_{\hrho_a}^{l-2}G(a,a)_k\Phi^0
\endalign$$
\qed\enddemo 
We determine $a_k$ such that 
$$
\frac1{k!}d\hrho_{a_k} \Phi^0 = - Ob_a( a_{<k})
\tag{$eq_k$}$$
In order that there exists a solution of $eq_k$, it is necessary that 
$[Ob_k] =0 \in H^2(\#)$.  If $[Ob_k] =0 $, we define $a_k$ by 
$$
\frac1{k!}\hrho_{a_k}\phi^0 = -d^*_1 G_\# (Ob_k(a_{<k} )) .
\tag2-2-8$$
Since $Ob_k(a_{<k}) $ is d-exact,
then we also have a criterion,
\proclaim{Theorem 2-2-2} 
If $p^2$ is injective, 
then $Ob_k (a_{<k})$ vanishes for all $k$.
\endproclaim
Thus we construct a power series $a(t)$ satisfying 
$ d\rho_{g(t)} \Phi^0=0$. 
Next we must prove that this power series $a(t)$ converges for sufficiently small $t$. 
\head \S2-3 the convergence
\endhead
We rewrite definition 1-6 by using the Sobolev norm.
\proclaim{Definition 1-6} 
A closed element $\Phi^0 \in \E_s(X)$ is unobstructed if there 
exists an integral curve $\Phi_t(a)$ in $\wtil{\M}_s(X)$ for each $a \in
\E_s\cap \Cal H$ such  that 
$$
\frac d{dt}\Phi_t(a)|_{t=0} = a
$$
An orbit $\O$ is unobstructed if any $\Phi^0 \in \wtil{\M}_s(X)$  is
unobstructed  for every compact $n$ dimensional manifold $X$
\endproclaim
The rest of this subsection is devoted to the proof theorem 1-7 ( criterion
of  unobstructedness). 
Our method is similar to the one of the Kodaira-Spencer theory. 
(See the extremely helpful book by Kodaira [18] for technical details.)
\demo{Proof of theorem 1-7}
We already have a formal power series $a(t)$ such that 
$$
d\rho_{g(t)} \Phi^0 =0.
$$ 
Hence it is sufficient to prove that $a(t)$ uniformly converges with
respect to the Sobolev norm $\| \,\|_s$.
Since $(L_a)_k\Phi^0 = L_{a_k}\Phi^0=-d\hrho_{a_k}\Phi^0$ and $dR_k=0$,
$a_k$ satisfies
$$
-\frac1{k!}d\hrho_{a_k}\Phi^0 = Ob_k.
\tag2-3-1$$
As in section 2-2, $Ob_k$ is an element of $\Gam(E^2)$.
$Ob_k$ is also written as
$$
Ob_k =\frac1{2!}d\hrho_{a_{\ss-style{<k} }}^2\Phi^0 +\cdots +\frac1{(k-1)!}d\hrho_{a_{\ss-style{<k} }}^{k-1}\Phi^0
\tag2-3-2$$
By (2-3-2), we see that $Ob_k$ is an exact form. 
Hence if the map $p^2 \: H^2( \#) \to \oplus_i H^{p_i+1}(X)$
is injective, then the class $[Ob_k]\in H^2(\#)$ vanishes.
Hence we obtain a solution of the equation (2-3-1) by 
$$
\frac1{k!}\hrho_{a_k}\Phi^0 = -d_1^* G_\# (Ob_k) \in E^1.
\tag2-3-3$$
We assume that $a_k $ belongs to the orthogonal complement of 
the Lie algebra of the isotropy group $H$.
Hence $a_k$ is defined uniquely by $\hrho_{a_k}\Phi^0$ and 
we have the estimate 
$$
\| a_k \|_s = C_1 \| \hrho_{a_k}\Phi^0\|_s
\tag2-3-4$$
Hence by (2-3-3), we have a formal power series,
$$
a=\sum_{k=1}^\infty \frac 1{k!} a_k t^k.
$$
Given two power series $P(t)=\sum_k p_k t^k$ and  $Q(t)= \sum_k q_k t^k$,  
if $p_k < q_k$ for all $k$, we denote it by 
$$
P(t) \ll Q(t).
$$
We denote by $(P)_k$ the homogeneous part of degree $k$ of $P(t)$.
Let $A(t)$ be a convergent series given by 
$$
A(t) = \frac{b}{16c}\sum_{k=1}^\infty\frac{c^k t^k}{k^2},
$$
with $b>0, c>0$. $b$ and $c$ will be determined later. 
As regards $A(t)$ we have the following inequality 
(see section 5-3 in [18]), 
$$
A(t)^l \ll(\frac{c}{b})^{l-1} A(t).
$$
Fix a natural number $s$. We shall show by induction on $k$ 
if we choose appropriate large $b$ and $c$, 
$$
\| a_{\leq k} \|_s \ll A(t),
\tag{ $*_k$ }
$$
where $\|a_{\leq k} \|_s = \sum_{l=1}^k \frac1{l!}\|a_l \|_s t^l$
We assume $*_{k-1}$ holds and make an estimate $\| a_k\|_s$.
By (2-3-3) we have the inequalities for constants $C_2, C_3$, 
$$\align
\frac1{k!}\| a_k \|_s = &C_1 \frac1{k!}\| \hrho_{a_k} \Phi^0 \|_s 
=C_1 \| d^*_1 G_\# (Ob_k) \|_s \\
<&  C_2 \| G_\# (Ob_k) \|_{s+1}
<C_3 \| Ob_k \|_{s-1}
\endalign$$
By theorem 2-2-1, we have an estimate, 
$$\align
&\s-style{
\| Ob_k \|_{s-1} <
\(\frac1{2!} \| G(a,a) \Phi^0\|_{s-1} 
+\frac1{3!} \| Ad_{\hrho_a}G(a,a)\Phi^0 \|_{s-1} + 
\cdots+ \frac1{k!}\|Ad_{\hrho_a}^{k-2}G(a,a)\Phi^0\|_{s-1}\)_k}\\
<&\s-style{
C_4\(\frac1{2!}\| G(a,a)\|_{s-1}+\frac1{3!}2 \|
a\|_{s-1}\|G(a,a)\|_{s-1}+\cdots\frac1{k!}2^{k-2}\|a \|_{s-1}^{k-2}
\|G(a,a)\|_{s-1}\)_k} \\ 
<&\s-style{C_4 ( \til{f}\( 
2\|a_{\ss-style<k}\|_{s-1}\)\,\|G(a,a)_{\ss-style k}\|_{s-1}})_k,
\endalign$$
where $\til{f}(x) =\frac 1{x^2}( e^x -1 -x)$.
We have an estimate of $G(a,a)$, 
$$
\|G(a,a)\|_{s-1} < C_5 \| a\|_s \|a\|_s
\tag 2-3-5$$
Hence by (2-3-5),
$$
\| Ob_k \|_{s-1} < C_6\(\big(\frac1{2!} 
+ \frac1{3!}2\|a_{\ss-style<k}\|_s+\cdots 
+\frac1{k!}2^{k-1}\|a_{\ss-style<k}\|_s^{k-2}
\big) \|a_{\ss-style<k}\|_s^2\)_k,
$$
where $C_6$ is a constant.
By the hypothesis of the induction, 
$$\align
&\s-style{
\| Ob_{ k} \|_{s-1} < C_6
\big( ( \frac1{2!} 
+ \frac1{3!}2A(t)+\cdots 
+\frac1{k!}2^{k-1}A(t)^{k-2} )A(t)^2 }\big)_k\\
&<\s-style{
C_6 \big( (  \frac1{2!} 
+ \frac1{3!}2A(t)+\frac1{4!}2^2(\frac{b}{c})A(t)\cdots 
+\frac1{k!}2^{k-1}(\frac{b}{c})^{k-3}A(t) )(\frac{b}{c})A(t)}\big)_k\\
&\s-style{
=C_6 \big( (\frac1{2!} (\frac{b}{c})A(t)
+ \frac1{3!}2(\frac{b}{c})A(t)+\frac1{4!}2^2(\frac{b}{c})^2
A(t)\cdots 
+\frac1{k!}2^{k-1}(\frac{b}{c})^{k-2}A(t))}\big)_k\\
&\s-style{
<C_6 \frac 1{2p}( e^{2p}-1 -2p ) A_k(t)},
\endalign$$
where $p =\frac{b}{c}$.
We define $p$ by  $C_6 \frac 1{2p}( e^{2p}-1 -2p )=1$. 
Then we obtain 
$$
\| Ob_{k} \|_{s-1} < A_k(t)
$$
Therefore we have 
$$
\frac1{k!}\| a_k \|_s < C_3 A_k(t).
$$
Since $A(t)$ is a convergent series for sufficiently small $t$, 
we see that $a(t)$ uniformly convergents.
\qed\enddemo
Further we assume that 
$$\align
&d\hrho_{a_1}\Phi^0 =0,\\
&d^*_0 \hrho_{a_1}\Phi^0 =0,
\endalign
$$
where $d^*_0$ is the adjoint operator and
$$
\CD 
0@>>>E^0 @> d_0>> E^1 @>d_1>>\cdots.
\endCD
$$
We also apply
the elliptic regularity to $\rho_{g(t)} \Phi^0$.  As in our
construction, we have 
$$\align
d\rho_{g(t)}\Phi^0 =&\hrho_{a}\Phi^0 + \sum_{l=2}^\infty\frac1{l!} d\hrho_a^l\Phi^0 =0\\
d_0^*\hrho_a\Phi^0=0
\endalign
$$
Hence $\hrho_a\Phi^0$ is a weak solution of an elliptic differential equation, 
$$
\trian_\#\hrho_a\Phi^0 + d_0^*(\sum_{l=2}^\infty\frac1{l!} d\hrho_a^l\Phi^0
) =0,
$$
where $\trian_\#$ is the Laplace operator of the complex $\#$.
Hence we obtain 
\proclaim{Theorem 2-3-1} 
If $p^2$ is injective, then there exists a smooth solution of the equation
(eq$_*$) for all  tangent $[\h{\rho}_a\Phi^0] \in \H^1(\#_{\Phi^0})$. 
i.e., 
There exists a smooth form $\rho_{\exp{a(t)}}\Phi^0 \in \wtil{\M}(X)$ such that 
$$
\(\rho_{\exp{a(t)}}\Phi^0\)' |_{t=0} = \h{\rho_a}\Phi^0
$$
\endproclaim

\head
\S 3 Proof of theorems
\endhead
In this section we assume that an orbit $\O \in \oplus_i \w^{p_i}$ is metrical, elliptic and topological.
We shall show that the moduli space $\M_{\ss-style{\O}}(X)$ is a manifold.
\subhead 
\S3-1 
\endsubhead 
In this subsection we explain preliminary results, which  are related to functional analysis on manifolds. 
Our discussion will heavily depend on 
[7] and [22] and we use the same notation as in section one and two.
Let $X$ be a smooth $n$ dimensional manifold and  $F$ a smooth fibre bundle over $X$. Then we have a Hilbert manifold 
$L^2_s(F)$ consisting of those sections of $F$ 
which in local coordinates are defined by functions 
square integrable up to order $s$, i.e., Sobolev space 
$L^2_s(X)$. We also define a Banach manifold $C^k(F)$ by  functions $C^k(X)$ and  for $s > k+ \frac n2$, we have a smooth inclusion $L^2_s(F) \subset C^k(F)$.
If we consider the case $F=X\times X$, 
we find that the sections of $F$ are exactly the maps 
of $X$ to $X$. Then $C^1(F)$ is the set of $C^1$ maps 
from $X$ to $X$ with the topology of uniform convergence up to the first derivative. 
We define $C^1$Diff$(X)$ by the $C^1$ diffeomorphisms:
$$
C^1\text{Diff}(X)= \{ \, f \in C^1(F)\, |\, f^{-1}\in 
C^1(F)\, \}.
$$
We denote by $C^1\Diff(X)$ the identity component of 
$C^1$Diff$(X)$. 
Pick $s+1 > \frac n2 +1$.
We define Diff$^{s+1}_0$(X) by  
$$
\text{Diff}^{s+1}_0(X) = C^1\text{Diff}_0(X) \cap L^2_{s+1} (F).
$$
Then Diff$^{s+1}_0(X)$ is the Hilbert manifold. 
In section 3 of [7] it is shown that $\Diff^{s+1}(X)$ is a topological group 
under the operation of composition of mappings.
Further we have the action  $A$ by using pull back:
$$
A:\E_s \times \text{Diff}^{s+1}_0(X) \arrow \E_s,
$$
where $\E_s= C^1(\A(X))\cap L^2_s (\oplus_i \w^{p_i})$ ( see  2-1-1 in section two).
Then the action $A$ is well defined and continuous .
For an element $\Phi\in \E_s$, we define 
$ A_\Phi\: \Diff^{s+1}(X) \arrow \E_s$ by 
$A_\Phi ( f) = f^*\Phi$. 
This map is continuous. Furthermore if $\Phi$ is smooth, then $A_\Phi$ is also smooth. 
We define the moduli space $\M(X)$ consisting of smooth forms  as in section one.
We shall extend this moduli in terms of Sobolev Space.
We define $\M_s(X)=\M_{\ss-style{\O},s}(X) $ by 
$$
\M_s(X)= \widetilde{\M}_s (X) /\text{Diff}_0^{s+1}(X),
$$where
$$
\widetilde{\M}_s (X) = 
\{ \,\Phi \in \E_s \, |\, d\Phi =0 \, \}.
$$
We denote by $\pi$  the natural projection 
$$
\pi\:\widetilde{\M}_s (X) \arrow \M_s(X).
$$
\subhead 
\S 3-2
\endsubhead
In section two, 
we construct the family of smooth closed forms
$\rho_{g(t)}\Phi^0$ parametrized by $t$, where  
$g(t) =\exp a(t)$ is written as
$$
\align
&a(t) = \sum_{k=1}^\infty \frac1{k!} a_k t^k,\\
&\hrho_{a_1}\Phi^0\in
\H^1(\#_{\Phi^0}).
\endalign
$$ 
Substitute $t=1$ into $\rho_{g(t)}\Phi^0$, we have 
a map $\til\kappa$, 
$$\align
\til{\kappa}\: & S_{\Phi^0} \to  \wtil{\M}(X),\\
\hrho_{a_1}&\Phi^0 \mapsto \rho_{g(1)}\Phi^0,
\endalign$$
where $S(=S_{\Phi^0})$ is a sufficiently small open set 
of $\H^1(\#_{\Phi^0})$.  Let $\pi$ be the natural projection $\wtil{\M}(X) \to \M(X)$. We define $\kappa$ as 
$$
\kappa=\pi \circ \til\kappa 
\: S \to \M(X).
$$
Then we shall show 
\proclaim{Theorem 3-2-1} 
We assume that $p^1\: H^1(\#_{\phi^0})\to 
\oplus_i H^{p_i}(X)$ is injective.
Then $\k\: S \to \M(X)$ 
is injective for a sufficiently small open set $S \subset 
\H^1(\#_{\phi^0})$, 
\endproclaim
\demo{Proof} 
We denote by $\Cal S(=\Cal S_{\Phi^0})$ the image $\til\k(S)$. We call $\Cal S$
a slice, which is a transversal submanifold for the action of $\Diff(X) $ on
$\wtil{\M}(X)$. Since the action of 
$\Diff(X)$ preserves the de Rham cohomology class, 
taking the cohomology class, 
we have the map 
$$
P\: \M(X) \arrow \oplus_i H^{p_i}_{dR}(X).
$$ 
The slice $\Cal S_{\Phi^0}$ is parametrized by an open set 
$S$. We denote by $P|_{\Cal S}$ the restricted map to 
$\Cal S$. 
Then the differential $dP|_{\Cal S}$ at $\Phi^0$ is given by $p^1\: H^1(\#_{\Phi^0}) \to \oplus_iH^{p_i}_{dR}(X)$. Since $p^1$ is injective, it follows that 
$P|_{\Cal S}$ is locally injective.  
Hence there exists an open set $S\in \H^1(\#)$ such that 
$P\circ \k \: S \arrow \oplus_i H^{p_i}_{dR}(X)$
is injective. We assume that there exists a diffeomorphism $f \in \Diff(X)$ such that 
$$
f^* \Phi^1 = \Phi^2,
$$
for some $\Phi^1, \Phi^2 \in S_{\Phi^0}$.
By taking each cohommology class, we have 
$$
P(f^*\Phi^1) = P(\Phi^2). 
$$
Since $\Diff(X)$ trivially acts on cohomology groups, 
we also have 
$$
P(f^*\Phi^1)=P(\Phi^1).
$$
Hence we have 
$$
P(\Phi^1)=P(\Phi^2).
$$
Since $P|_{\Cal S}$ is injective, we see that $\Phi^1=\Phi^2$. Hence $\k$ is injective.
\qed\enddemo
Next we shall show that $\k$ is surjective:
\proclaim{Theorem 3-2-2} 
If we take a sufficiently small open set  $U_{\Phi^0}$ of 
$\pi(\Phi^0)$ in $\M(X)$, 
the map $\kappa\: S \to U_{\Phi^0}$ is surjective 
for an open set $S$ in $\H^1(\#)$.
\endproclaim
In order to prove theorem 3-2-2, we must explain 
the completion with respect to Sobolev norm and the following lemma 3-2-3 and proposition 3-2-4.
Let $\xi$ be a vector field on $X$. 
We assume that $\xi \in C^1(TX)\cap L^2_{s+1}(TX)$. 
Then there is the diffeomorphism $f_\xi$ corresponding to $\xi$, 
where $f_\xi \in \Diff^{s+1}(X)$ (see 3-1).
Since $\E_s(X)$ is invariant under the action of diffeomorphism, we have the action of $f_\xi$ on $\E_s(X)$, i.e., 
$\rho_{e^a}\Phi^0 \to f_\xi^* \rho_{e^a}\Phi^0,$ 
where $a\in $End$(TX)$ and $\rho_{e^a}\Phi^0\in 
\E_s(X)$. Hence there exists $b_\xi(=b_{f_\xi})\in $End$(TX)$ such 
that 
$$
f_\xi ^* \rho_{e^a}\Phi^0 = \rho_{\exp b_{\ss-style{\xi}}}\Phi^0.
\tag1$$
We assume that $\rho_{e^a}\Phi^0$ is smooth with $d\rho_{e^a}\Phi^0 =0$, i.e., $\rho_{e^a}\Phi^0 
\in \wtil{\M}(X)$. Then we 
shall show that for a sufficiently small 
$a$, there exists a vector field $\xi$ satisfying 
$$
d_0^*\hrho_{b_{{\xi}}}=0
\tag2$$
Since $\hrho$ is the differential representation of $\rho$, 
$\hrho_{{b_{\xi}}}\Phi^0$ is written as 
$$
\hrho_{{b_\xi}}\Phi^0 = \rho_{\exp b_\xi}\Phi^0 -
\Phi^0 - \sum_{k\geq 2}\frac1{k!}\hrho_{b_\xi}^k 
\Phi^0.
\tag3$$
Since $\rho_{\exp b_\xi}\Phi^0= f_\xi^*\rho_{e^a}\Phi^0$ 
and $f^*\rho_{e^a}\Phi^0$ is closed, 
by using the identity: $L_\xi = d\circ i_\xi + i_\xi \circ d$,  we have 
$$\align
\rho_{\exp b_{\xi}}\Phi^0 =
&f_\xi^* \rho_{e^a}\Phi^0=
\rho_{e^a}\Phi^0+di_\xi \rho_{e^a}\Phi^0 + \cdots,\tag4\\
=&\Phi^0+\hrho_a \Phi^0 +di_\xi \Phi^0 + H( \xi ,a ),
\endalign$$where 
$H(\xi,b )$ denotes the higher order terms with respect to $\xi $ and $a$.
Substituting (4) into (3), we have 
$$
\hrho_{b_{\xi}}\Phi^0 = \hrho_a \Phi^0 + di_\xi \Phi^0 +W( \xi, a),
\tag5$$where $W(\xi, a)$ denotes the higher order terms.
In order to solve the equation (2), we need the following estimate of $W(\xi,
a)$ for sufficiently small 
$\xi_1, \xi_2$ and $a$:
\proclaim{Lemma 3-2-3} 
$$
\align &\| W( \xi_1,a ) - W(\xi_2, a) \|_{s}  < \e \| \xi_1- \xi_2 \|_{s+1},
\endalign 
$$
where $\e< 1$ is a constant.
\endproclaim
\demo{Proof} 
We have the action of the diffeomorphism Diff$_0^s(X)$ by 
$$
\Diff^{s+1}(X)\times L_s^2(\oplus_i\w^{p_i}) 
\arrow L_{s}^2 (\oplus_i\w^{p_i} ).
\tag6$$ 
We denote this map by $A$. 
Then we have the map $A_\Phi\: \Diff^{s+1}(X) \to \E_s$
by $A_\Phi (f) = f^*\Phi$ as in section 3-1. 
Then $A_\Phi$ is smooth if $\Phi$ is smooth. 
From our assumption $\rho_{e^a}\Phi^0 \in 
\wtil{\M}(X)$, $a$ is smooth. 
We identify $\xi$ with $f_\xi \in \Diff(X)$.
Hence we have 
$$
\|A(\xi_1,\rho_{e^a}\Phi^0) - A(\xi_2, \rho_{e^a}\Phi^0) \|_{s} 
\leq \| dA \|_0 \| \xi_1 - \xi_2 \|_{s+1},
\tag7$$
where $dA$ denotes the differential of $A$, 
$\|dA \|_0$ is the $C^0-$ norm of $dA$.
Since $W$ is essentially written in terms of $A$,
from (7) we also have 
$$
\| W(\xi_1,a) - W(\xi_2, a) \|_{s} 
\leq \| dW \|_0 \| \xi_1 -\xi_2 \|_{s+1}.
\tag8$$
Since the differential of $W$ at the $(\xi,b)=(0,0)$ vanishes, 
we have 
$$
\|W(\xi_1, a) - W(\xi_2, a) \|_{s} 
\leq \e \| \xi_1 - \xi_2 \|_{s+1},
\tag9$$
where $\e<1$ is a constant.
Hence  by (9), we have the result.
\enddemo
\proclaim{ Proposition 3-2-4} 
There exists a sufficiently small $\varepsilon>0$ satisfying the following: \par 
 For any $\rho_{e^a} \Phi^0 \in \wtil{\M}(X)$ 
 with  $\| a\|_s <\e$, 
there exists a smooth vector field $\xi $ such that 
$d_0^* \hrho_{b_{\xi}} \Phi^0=0$. 
\endproclaim
\demo{Proof} 
Let $\xi$ be a vector field and $f_\xi$ the diffeomorphism corresponding to $\xi$, 
where $\xi\in C^1(TX)\cap L^2_{s+1}(TX)$ and 
$f_\xi \in  \Diff^{s+1}(X)$ as before.
We shall construct a vector field $\xi$ satisfying 
the equation (2). 
By (5), the equation (2) is written as 
$$
d_0^* d i_\xi \Phi^0 + d_0^* \hrho_a\Phi^0 
+ d_0^*W(\xi, a )=0.
\tag10$$
We recall the complex $\#_{\Phi^0}$:
$$
\CD 
0@>>>E^0_{\Phi_0} @>d_0>>E^1_{\Phi_0} @>d_1>> 
E^2_{\Phi_0} @>d_2 >>\cdots, 
\endCD
\tag$\#_{\Phi_0}$
$$
where $d_0^*$ denotes the adjoint operator of $d_0$,
$i_\xi \Phi^0 \in E^0$ and $\hrho_a \Phi^0 \in E^1$. Then by using the Hodge
decomposition of $E^0$,
 we take $\xi$ such that 
the harmonic component of $i_\xi \Phi^0$ 
vanishes with respect to the complex $\#_{\Phi^0}$. 
Then the equation (10) is equivalent to the following:
$$
i_\xi \Phi^0+ G_\# d_0^* \hrho_a\Phi^0
+ G_\# d_0^*W(\xi, a) =0 
\tag11
$$ 
where $G_\#$ denotes the Green operator with respect to the complex $\#_{\Phi^0}$.
Then  given $a$ and $\Phi^0$, we shall show that there exists a solution $\xi$ of
(11). We denote by  Ker$\,\Phi^0$ the subbundle of TX
given by 
$$
\text{Ker} \,\Phi^0 =\{ \, \xi \in TX \, |\, 
i_\xi\Phi^0 =0 \, \}.
$$
Let Ker$^\perp$ be the orthogonal complement of 
Ker$\,\Phi^0$ in $TX$.
At first we define $\xi_1\in $Ker$^\perp$ by 
$$i_{\xi_1}\Phi^0=-d_0^*G_\#\hrho_a \Phi^0.
\tag12$$ 
Note that since the image of $d_0^*$ is in $E^0$, 
there is a unique vector field $\xi_1$ satisfying $(12)$.
Secondly we define $\xi_2\in$Ker$^\perp$ by 
$$
i_{\xi_2}\Phi^0= -G_\#d_0 \hrho_a\Phi^0 - G_\#d_0^* W( \xi_1, a ) \Phi^0.
\tag13$$
Inductively we define $\xi_k\in$Ker$^\perp$ by 
$$
i_{\xi_k}\Phi^0 = - G_\# d_0 \hrho_a\Phi^0 - G_\# d_0^* W( \xi_{k-1}, a ) \Phi^0.
\tag14$$
Since $\xi_k \in$Ker$^\perp$, we have an estimate 
$$
\| \xi_k \|_{s+1} = C\|i_{\xi_k}\Phi^0\|_{s+1},
\tag15$$
where $C$ denotes a constant.
Then by (14) and  the elliptic estimate of $G_\#$ and $d_0^*$, we have
$$\align
\| \xi_{k+1} - \xi _k \|_{s+1} =&  
\| \big(G_\#d_0^*  W(\xi_k, a) - G_\#d_0^* W(\xi_{k-1},a)\big)\Phi^0\|_{s+1} 
\tag16\\
\leq&C_1 \| \big(W(\xi_k, a) - W(\xi_{k-1}, a ) \big)\|_{s} \\
\endalign
$$
where $C_1$ is a constant.
Hence by lemma 3-2-3, we have 
$$
\| \xi_{k+1} - \xi_{k}\|_{s+1} = \e \| \xi_k - \xi_{k-1}\|_{s+1},
\tag17$$
where $\e<1$ is a constant.
 Hence it follows from (17) that 
 the sequence $\{ \xi_k \}$ uniformly converges to 
 some $\xi_\infty$ with respect to the norm $\|\, \|_{s+1}$. Then by (14) we see that $\xi_\infty$ is 
 a solution of the equation (11).
 Hence we have a vector field satisfying (2).
 \qed\enddemo
\demo{Proof of theorem 3-2-2} 
By proposition 3-2-4, if we define an  open set $U_{\Phi^0}$ by  the image: 
$$
U_\e({\Phi^0})=\pi\(\{ \, \rho_{e^a}\Phi^0 
\in \wtil{\M}(X)\, |\,\, \|a \|_s < \e , \, \,\}\),
$$ 
then there exists a diffeomorphism $f_\xi$ such that 
$$\align
&f_\xi^* \rho_{e^a}\phi^0 = \rho_{\exp b_\xi}\Phi^0\\
&d_0^* \rho_{b_\xi}\Phi^0 =0.
\endalign 
$$
Hence $\rho_{\exp b_\xi}\Phi^0$ is in the image of 
$ \til{\kappa}(S)$. Hence it follows that $\kappa$ is surjective.
\qed\enddemo
\subhead 
\S 3-3
\endsubhead
Let $\O$ be an orbit in $\oplus_i\w^{p_i}$ 
as in section one. Then we have the moduli space 
$\M_{\ss-style\O}(X)(=\M(X))$ as the quotient space 
$\wtil{\M}(X)/\Diff(X)$. 
\proclaim{Proposition 3-3-1} 
We assume that the orbit $\O$ is elliptic, metrical and 
topological.
Then the quotient $\M_{\ss-style{\O}}(X)$ is Hausdorff.
\endproclaim
\demo{Proof} 
Since the orbit $\O$ is metrical, we have the metric $g_\Phi$ for every $\Phi \in \E_s(X)$. 
Then each tangent space 
$T_\Phi \E_s(X)\subset 
L^2_s( \oplus_i \w^{p_i})$ has the $L^2$ metric 
in terms of $g_\Phi$. Hence it gives $\E_s(X)$
a smooth Riemannian structure
(see section 4 in [7]). 
The fundamental property of the Riemannian structure 
on $\E_s(X)$ is that it is invariant under the action of $\Diff^{s+1}(X)$, so that is, 
$\Diff^{s+1}(X)$ acts on $\E_s(X)$ isometrically.  
Since $\wtil{\M}_s(X)$ is the intersection 
$\Cal H\cap \E_s(X)$, we have the induced distance 
on $\wtil{\M}_s(X)$ from the Riemannian structure on 
$\E_s(X)$. 
We denote by $d$ the distance on
$\wtil{\M}_s(X)$ and $\pi$  the natural projection
$\pi\: \widetilde{\M}_s(X)\to \M_{s}(X)$. Then we define $ d(\pi(\Phi^1),\pi(\Phi^2) )$ by 
$$
d(\pi(\Phi^1) ,\pi(\Phi^2)) = \inf_{\ss-style{f,g \in \Diff^{s+1}(X)}} d ( f^*\Phi^1,
\,\,g^*\Phi^2 ),
\tag3-3-1$$
where $\Phi^1,\Phi^2\in \wtil{\M}_s(X)$.
For simplicity we denote by $\Cal D$ 
the group $\Diff^{s+1}(X)$.
Since the action of $\Diff^{s+1}(X)(=\Cal D)$ preserves the distance $d$, we have 
$$
d( \pi(\Phi^1),\pi(\Phi^2) ) = 
\inf_{\ss-style{\,\,f \in
\Cal D }}d ( f^*\Phi^1, \,\,\Phi^2).
\tag3-3-2$$
Hence we have the triangle inequality , 
$$
\align
&\quad\quad  d(\pi (\Phi^1), \,\,\pi(\Phi^2) ) + d( \pi(\Phi^2),\,\,\pi(\Phi^3)) 
\tag3-3-3\\
=&\inf_{\ss-style{f \in \Cal D}}d( f^*\Phi^1,\Phi^2) + 
\inf_{\ss-style{g\in \Cal D} }d(\Phi^2,
g^*\Phi^3)\\ 
\leqq &\inf_{\ss-style{f,g\in \Cal D}} d( f^*\Phi^1,\,\,g^*\Phi^3) = d(
\pi(\Phi^1),\,\pi(\Phi^3)).
\endalign
$$
We shall show that $d$ induces a distance of $\M(X)$. 
We assume that the  $d(\,\pi(\Phi^0),\,\,\pi(\Phi)\,)=0$ for smooth elements $\Phi, \Phi^0 \in 
\wtil{\M}(X)$.
Then by (3-3-1), we have 
$$
{\inf_{\ss-style{f\in\Diff(X)}}d(\Phi^0, f^*\Phi)=0}.
\tag3-3-4
$$ Hence $f^*\Phi$ is in a small
neighborhood 
$U_\e(\Phi^0)$ at $\Phi^0$.  By  theorem 3-2-3 and proposition 3-2-4,  
there exists a diffeomorphism $f_\xi$ such that 
$$
f_\xi^* (f^*\Phi^0) \in \Cal S_{\Phi^0},
\tag3-3-5$$
where $\Cal S_{\Phi^0}$ is the family parametrized by 
an open set of harmonic forms $\H^1(\#_{\Phi^0})$.
 We define the distance $d_{\H^1(\#)}$ on $\H^1(\#_{\Phi^0})$ by using harmonic forms in terms of $g_{\Phi^0}$.  There is the distance $d_{H_{dR}}$ on  the direct sum of the de Rham cohomology groups $\oplus_i H^{p_i}(X)$
 by  using harmonic representations with respect to $g_{\Phi^0}$. Since $p^1$ is injective, $( H^1(\#_{\Phi^0}),
 d_{\H^1(\#) } )$ is isometrically embedded into 
 $(\oplus_i H^{p_i}(X), d_{\H_{dR}} )$.  
Since $p^1$ is injective, we have the injective map 
$$
P|_\Cal S\: S\to \oplus_i H^{p_i}_{\ss-style{dR}}(X).
\tag3-3-6$$
Since the differential $dP|_{\Cal S}\: T_{\Phi^0}\Cal S (=\H^1(\#_{\Phi^0}) )\arrow \oplus_i H^{p_i}(X)$ 
is isometric,  we have that 
$$
d( \Phi^0 , f_\xi^*( f^*\Phi) ) > 
C d_{\ss-style{dR}}
\(P(\Phi^0), P( f_\xi^*(f^*\Phi) )\),
\tag3-3-7$$
where $C$ is a positive constant.
(Note that the distance $d$ restricted to the slice 
$\Cal S$ is (locally) equivalent to the distance $d_{\H^1(\#)}$.)
Since $\Diff(X)$ acts on $\oplus_i H^{p_i}(X)$ trivially, 
we have
$$
d( \Phi^0 , f_\xi^*( f^*\Phi) ) > 
C d_{\ss-style{dR}}
\(P(\Phi^0), P( \Phi) )\),
\tag3-3-8
$$
where $C$ does not depend on $f$ and $\Phi$.
Hence
$$
\inf\Sb {\ss-style{f\in\Diff(X)}}\endSb d(\Phi^0, f^*\Phi)>
 C\,\, d_{\ss-style{dR}}(
P(\Phi^0),\,\,P(\Phi)).
\tag3-3-9$$
 Hence from our assumption (3-3-4), we have $P(\Phi) = P(\Phi^0)$. 
Since $P|_\Cal S$ is injective, 
$\Phi^0 = f_\xi^*( f^* \Phi)$. Hence we have
$\pi(\Phi)=\pi(\Phi^0)$. Hence $d$ is a distance on
$\M(X)$.  
\qed\enddemo
\subhead
\S 3-4 Proof of main theorems
\endsubhead
\demo{Proof of theorem 1-8}
Since the orbit $\O$ is elliptic and topological, 
we have the slice $S_{\Phi^0}$ as 
coordinates of the moduli space $\M(X)$
by  theorem 3-2-1 and 3-2-2 in section 3-2.
Hence $\M(X)$ is a manifold.
Since $\O$ is metrical, $\M(X)$ is Hausdorff 
by proposition 3-3-1 in section 3-3.
The slice $S_{\Phi^0}$ is homeomorphic to an open set of the cohomology group H$^1 (\#_{\Phi^0} )$. 
\qed\enddemo
\demo{Proof of theorem 1-9}
The  slice $S_{\Phi^0}$ is local coordinates of $\M_{\ss-style{\O}}(X)$. 
The differential $dP$ coincides with the injective map $p^1$. Hence $P$ is locally injective.
\qed\enddemo
Since $\O$ is metrical, we have the metric $g_{\Phi}$ for each $\Phi \in \E$.
Hence the metric $g_{\Phi}$ defines the metric on each tangent space $E^1 =
T_\Phi \E$.  Then $\E$ can be considered as a Riemannian manifold.
Then we see that the action of Diff$_0(X)$ on  $\E$ is isometry. 
Let $I({\Phi})$ be the isotropy group of Diff$_0(X)$ at $\Phi$, 
$$
I({\Phi}) = \{\, f \in \text{Diff}_0(X) \, |\, f^*\Phi =\Phi\, \}.
$$ 
Let $S_{\Phi^0}$ be a slice at $\Phi^0$. Then we shall compare $I({\Phi^0})$ to 
other isotropy group $I({\Phi})$ for $\Phi \in S_{\Phi^0}$.
\proclaim{Theorem 1-10}
Let $I({\Phi}^0)$ be the isotropy group of $\Diff(X)$ at $\Phi^0$ and $S_{\Phi^0}$ the slice at
$\Phi^0$.  Then $I({\Phi^0})$ is a  subgroup of the isotropy group $I({\Phi})$ for each $\Phi\in
S_{\Phi^0}$. ( We take $S_{\Phi^0}$ sufficiently small for necessary.)
\endproclaim
\demo{Proof of theorem 1-10}
From definition of $S_{\Phi^0}$, the slice $S_{\Phi^0}$ is invariant under the action of 
$I({\Phi^0})$. 
The restricted map $P|_{S_{\Phi^0}} \: S_{\Phi ^0}\to \underset i \to \oplus $H$^{p_i} (X)$ is locally injective. Since the action of Diff$_0(X)$ preserves each class of de Rham cohomology group, 
we see that the action of $I({\Phi^0})$ is trivial on the slice $S_{\Phi^0}$ 
for sufficiently small 
$S_{\Phi^0}$.
Hence $I({\Phi^0})$ is a subgroup of the isotropy group $I({\Phi})$ for each 
$\Phi \in S_{\Phi^0}$.
\qed\enddemo
\demo{Proof of theorem 1-11}
The slice $S_{\Phi^0}$ is local coordinates of $\M_{\ss-style{\O}}(X)$ and 
the action of $\Diff(X)$ on $\E$ is isometry.
Hence the moduli space $\wtil{\M}_{\ss-style{\O}}(X)/\Diff(X)$ is locally homeomorphic to 
the quotient space $S_{\Phi^0}/I(\Phi^0)$, where $I(\Phi^0)$ is the isotropy. 
Hence we see that there is an open set $V$ of $T_{\ss-style{\Phi^0}} S_{\ss-style{\Phi^0}}$
with the action of
$I(\Phi^0)$  such that the quotient $V/I(\Phi^0)$ is homeomorphic to $S_{\Phi^0}/I(\Phi^0)$.  
$T_{\Phi^0}S_{\Phi^0}$ is
isomorphic to $H^1(\#_{\ss-style {\Phi^0}})$ and the action of
$I(\Phi^0)$ on $H^1(\#_{\ss-style {\Phi^0}})$ is a isometry with respect to $g_{\ss-style {\Phi^0}}$. The
action of $I(\Phi^0)$ preserves the integral  cohomology class. It implies that the image of $I(\Phi^0)$ is a subgroup of
O$\,(H^1(\#_{\Phi^0}))\cap $ End$\,(\oplus_i ( H^{p_i}(X, \Z))$, 
where O$\,(H^1(\#_{\Phi^0}))$ denotes the orthogonal group.
Then we see that  $V/I(\Phi^0)$ is the quotient by a finite group. 
\qed\enddemo

\head \S4. Calabi-Yau structures
\endhead
\subhead 
\S4-1.SL$_n(\C)$ structures 
\endsubhead
Let $V$ be a real $2n$ dimensional vector space. 
We consider the complex vector space $V\otimes \C$ and 
a complex form $\Ome \in \w^n V^*\otimes\C$. 
The vector space ker\,$\Ome$ is defined as 
$$
Ker\, \Ome = \{ \, v \in V\otimes \C \, |\, i_v \Ome =0 \, \},
$$
where $i_v$ denotes the interior product.
\proclaim{Definition 4-1-1 (SL$_n(\C)$ structures) } 
A complex $n$ form $\Ome$ is a SL$_n(\C)$ structure on $V$ if 
$\dim_{\ss-style\C} $Ker\, $\Ome =n$ and Ker\, $\Ome \cap \ol{Ker\,\Ome}= \{0\}$,
where $\ol{Ker\, \Ome}$ is the conjugate vector space.
\endproclaim
We denote by $\A_{\ss-style{SL}}(V)$ the set of SL$_n(\C)$ structures on $V$. 
We define the almost complex structure $I_\Ome$ on $V$ by 
$$
I_\Ome (v) = 
\cases
-\sqrt{-1}v&\quad \text{ if  } v\in Ker\, \Ome,\\
\sqrt{-1}v&\quad \text{ if }v \in \ol{Ker\, \Ome}.
\endcases
$$
So that is, Ker\, $\Ome= T^{0,1}V$ and $\ol{Ker \Ome}=T^{1,0}V$ and
$\Ome$ is a non-zero $(n,0)$ form on $V$ with respect to $I_\Ome$.
Let $ \Cal J$ be the set of almost complex structures on $V$. Then 
$\A_{\ss-style{SL}}(V)$ is the $\C^*-$bundle over$\Cal J$. 
We denote by $\rho$ 
the action of the real general linear group $G=GL (V)\cong GL( 2n ,\R)$ on the
complex
$n$ forms, 
$$
\rho \: \text{GL}(V) \arrow \text{End}\,(\w^n (V\otimes\C)^* ).
$$
For simplicity we denote by $\w^n_\C$ complex $n$ forms.
Since $G$ is a real group, $\A_{\ss-style{SL}}(V)$ is invariant under the action of $G$. 
Then we see that the action of $G$ on $\A_{\ss-style{SL}}(V)$ is transitive. 
The isotropy group  $H$ is defined as 
$$
H =\{\, g \in G\, |\, \rho_g \Ome = \Ome \, \}.
$$
Then we see $H=$SL$(n, \C)$. Hence the set of SL$_n(\C)$ structures $\A_{SL
}(V)$ is the homogeneous space, 
$$
\A_{\ss-style{SL}}(V) = G/H = GL(2n, \R) /SL(n,\C).
$$
(Note that the set of almost complex structures $\Cal J=$GL$(2n,\R)/$GL$(n,\C)$. )
An almost complex structure $I$ defines a complex subspace $T^{1,0}$ of dimension
$n$. Hence we have the map $\Cal J \arrow $Gr$(n, \C^{2n})$. We also have the
map  from $\A_{\ss-style{SL}}(V)$ to the tautological line bundle $L$ over the
Grassmannian  Gr$(n,\C^{2n})$ removed $0-$section.
Then we have the diagram: 
$$
\CD 
\A_{\ss-style{SL}}(V) @>>> L\backslash{0}\\
@VV\C^*V @VVV\\
\Cal J @>>>Gr(n,\C^{2n} )
\endCD
$$
$\A_{\ss-style{SL}}(V)$ is embedded as a smooth submanifold in $n-$forms $\w^n$. 
This is Pl\"ucker embedding described as follows, 
$$
\CD
\A_{\ss-style{SL}}(V) @>>> L\backslash{0}@>>>\w^n\backslash\{0\}\\
@VV\C^*V @VVV@VVV\\
\Cal J @>>>Gr(n,\C^{2n} )@>>>\CP^n.
\endCD
$$
Hence the orbit $\O_{\ss-style{SL}}=\A_{\ss-style{SL}}(V)$ is a submanifold in $\w^n$ 
defined by Pl\"ucker relations.
Let $X$ be a real $2n$ dimensional compact manifold.
Then we have the $G/H$ bundle $\A_{\ss-style{SL}}(X)$ over $X$ as in section 1. 
We denote by $\E=\E^1_{\ss-style{SL}}$ the set of smooth global sections of
$\A_{\ss-style{SL}}(X)$.  Then we have the almost complex structure $I_\Ome$
corresponding to 
$\Ome\in \E^1$. Then we have 
\proclaim{Lemma 4-1-2}
If $\Ome \in \E^1$ is closed, then the almost complex structure $I_\Ome$ is 
integrable.
\endproclaim
\demo{Proof}
Let $\{\theta_i\}_{i=1}^n$ be a local basis of $\Gam(\w^{1,0})$ with respect to 
$\Ome$. From Newlander-Nirenberg's theorem it is sufficient to show that 
$d\theta_i \in \Gam (\w^{2,0}\oplus \w^{1,1} )$ for each $\theta_i$. 
Since $\Ome$ is of type $\w^{n,0}$, 
$$
\theta_i \w \Ome =0.
$$
Since $d\Ome=0$, we have 
$$
d\theta_i \w \Ome =0.
$$
Hence $d\theta_i \in \Gam(\w^{2,0}\oplus \w^{1,1} )$.
\qed\enddemo
Then we define the moduli space of SL$_n(\C)$ structures on $X$ by
$$
\M_{\ss-style{SL}}(X) = \{\, \Ome \in \E^1_{\ss-style{SL}}\, |\, d\Ome =0 \, \}/
\text{Diff}_0(X).
$$ 
From lemma 4-1-2 we see that $\M_{\ss-style{SL}}(X)$ is the $\C^*-$bundle over 
the moduli space of integrable complex structures on $X$ with trivial canonical
line bundles.
\proclaim{Proposition 4-1-3}
The orbit $\O_{\ss-style{SL}}$ is elliptic. 
\endproclaim
\demo{Proof}
Let $\w^{p,q}$ be $(p,q)-$forms on $V$ with respect to  $I_{\Ome^0} \in
\A_{\ss-style{SL}}(V)$.  In this case we see that 
$$\align 
E^0 &= \w^{n-1,0}\\
E^1 &=\w^{n,0}\oplus \w^{n-1,1}\\
E^2 &=\w^{n,1}\oplus \w^{n-1,2}.
\endalign
$$
Hence we have the complex :
$$\CD
 \w^{n-1,0} @>\w u >>\w^{n,0}\oplus \w^{n-1,1}@>\w u >> 
\w^{n,1}\oplus\w^{n-1,2}@>\w u>>\cdots,
\endCD
$$
for $ u \in V$. 
Since the Dolbeault complex is elliptic, we see that 
the complex $0\arrow E^1 \arrow E^2 \arrow \cdots $
is exact. 
\qed\enddemo
\proclaim{Proposition 4-1-4}
Let $I_\Ome$ be the complex structure corresponding to $\Ome \in \E$. 
If $\partial \ol{\partial}$ lemma holds for the complex manifold $(X, I_\Ome)$,
then 

$$\align
&H^0(\#) \cong H^{n-1,0}(X),\quad H^1(\#) \cong H^{n,0}(X)\oplus
H^{n-1,1}(X),\\
&H^2(\#)\cong H^{n,1}(X)\oplus H^{n-1,2}(X) 
\endalign$$
and
$p^1,p^2$ are respectively injective, 
$$
p^1\:H^1(\#)\to H^n(X,\C),\quad
p^2\: H^2 (\#) \to H^{n+1}(X, \C)
$$ 
In particular, if
$(X, I_\Ome)$ is K\"ahlerian,
$p^k$ is injective for $k=1,2$.
\endproclaim
\demo{Proof}
As in proof of proposition 4-1-3 the complex \#$_\Ome$ is given as
$$
\CD 
\Gam( \w^{n-1,0} ) @>d>>\Gam ( \w^{n,0}\oplus\w^{n-1,1} ) @>d>>\Gam (\w^{n,1}
\oplus\w^{n-1,2} )@>d>>\cdots.
\endCD
$$
Then we have the following double complex:

$$\CD 
\Gam(\w^{n,0}
)@>\ol{\pa}>>\Gam(\w^{n,1})@>\ol{\pa}>>\Gam(\w^{n,2})@>\ol{\pa}>>\cdots\\
@A\pa AA @ A\pa AA@A\pa AA @.\\
\Gam(\w^{n-1,0})@>\ol{\pa}>>\Gam(\w^{n-1,1})@>\ol{\pa}>>
\Gam(\w^{n-1,2})@>\ol{\pa}>>\cdots\\
@A\pa AA @A\pa AA @A\pa AA @. \\
\Gam(\w^{n-2,0})@>\ol{\pa}>>\Gam(\w^{n-2,1})@>\ol{\pa}>>
\Gam(\w^{n-2,2})@>\ol{\pa}>>\cdots 
\endCD$$
Let $a = x+y$ be a closed element of $\Gam(\w^{n,1})\oplus\Gam(\w^{n-1,2})$. 
Then we have the following equations, 
$$\align 
&\ol{\partial}y =0 ,\tag 1\\
&\ol{\partial}x +\partial y =0. \tag 2 
\endalign
$$
Using the Hodge decomposition, we have 
$$
y= Har(y) +\ol{\pa}(\ol{\pa}^* G_{\ol{\pa}}\,y ),
\tag 3$$
where $G_{\ol{\pa}}$ is the Green operator with respect to the
$\ol{\pa}-$Laplacian
 and $Har (y)$ denotes the harmonic component of $y$. 
We also have 
$$
x = Har (x)  + \pa ( \pa^* G_\pa \,x ),
\tag 4$$
where $G_\pa$ is the Green operator with respect to the $\pa-$Laplacian 
and $Har (x)$ denotes the harmonic component of $x$. 
We put  $s =\pa^* G_\pa x$ and $t=\ol{\pa}^* G_{\ol{\pa}}y$ respectively.
Then we have from (2)
$$
\ol{\pa}\pa s+\pa \ol{\pa}t= 
\ol{\pa}\pa ( s-t) =0.\tag 5
$$
Applying $\pa\ol{\pa}$-lemma,   we
see from (5) that there exists a $\gam \in \w^{n-1,0}$ such that
$$
\pa (s-t) =\pa\ol{\pa}\gam.
\tag6
$$
Hence we have from (4), 
$$\align
x =& Har (x) + \pa s = Har (x) + \pa t +\ol{\pa}(-\pa\gam) \\
y =& Har (y) + \ol{\pa}t.
\endalign
$$
Thus if $Har(x)=0 $ and $ Har(y) =0$, then $a$ is written as 
$ a=x+y= d( t-\ol{\pa}\gam) $ where  $t-\ol{\pa}\gam \in 
E^1\cong \w^{n.0}\oplus\w^{n-1,1}$.  It implies that the map $p^2\: H^2 (\#)
\arrow H^{n+1} (X,\C)$ is injective and $H^2(\#) \cong H^{n,1}(X)\oplus
H^{n-1,2}(X)$.
\qed\enddemo
Hence from section 2, we have the smooth deformation space of SL$_n(\C)$ 
structures. However $\O_{\ss-style{SL}}$ is not metrical and the moduli space 
$\M_{\ss-style{SL}}(X)$ is not Hausdorff in general. 
In fact it is known that it is not Hausdorff for $K3$ surface. Hence in oder to obtain 
a Hausdorff moduli space, we must introduce extra geometric structures. 
The most natural structure is a Calabi-Yau structure.
\subhead \S4-2. Calabi-Yau structures
\endsubhead
Let $V$ be a real vector space of $2n$ dimensional. 
We consider a pair $\Phi=( \Ome, \ome)$ of a SL$_n(\C)$ structure $\Ome$ and 
a real symplectic structure $\ome$ on $V$, 
$$\align
\Ome &\in \A_{\ss-style{SL}}(V), \\
 \ome &\in \w^2 V^*, \quad\overset n \to{\overbrace 
{\ome\w
\cdots \w \ome} }\neq 0.
\endalign$$
We define $g_{\Ome,\ome}$ by 
$$
g_{\Ome,\ome}(u,v)=\ome(I_\Ome u, v), 
$$
for $u,v \in V$.
\proclaim{Definition 4-2-1(Calabi-Yau structures )} 
A Calabi-Yau structure on $V$ is a pair $\Phi=( \Ome,\ome) $ such that 
$$\align
&\Ome \w \ome =0 , \quad \ol{\Ome}\w \ome =0
\tag 1\\
&\Ome \w \ol{\Ome} = 
c_n \,\overset n\to{\overbrace{\ome\w \cdots \w\ome} 
}
\tag 2\\
&g_{\Ome,\ome}\text{ is positive definite.}\tag3
\endalign$$
where $c_n$ is a constant depending only on $n$,.i.e, 
$$
c_n =(-1)^{\frac{n(n-1)}2}\frac{2^n}{i^n n!}.
$$
\endproclaim
From the equation (1) we see that $\ome$ is of type $\w^{1,1}$ with respect to 
the almost complex structure $I_\Ome$. The equation (2) is called
Monge-Amp$\grave{e}$re equation. 
\proclaim{Lemma 4-2-2}
Let $\A_{\ss-style{CY}}(V)$ be the set of Calabi-Yau structures on $V$. 
Then There is the transitive action of $G=$GL$( 2n, \R)$ on
$\A_{\ss-style{CY}}(V)$  and $\A_{\ss-style{CY}}(V)$ is the homogeneous space 
$$
\A_{\ss-style{CY}}(V)= GL(2n,\R) / SU(n).
$$ 
\endproclaim
\demo{Proof}Let $g_{\Ome,\ome}$ be the K\"ahler metric. Then we have a unitary
basis  of $TX$. Then the result follows from (1) and (2).
\qed\enddemo
Hence the set of Calabi-Yau structures on $V$ is the orbit $\O_{\ss-style{CY}}$, 
$$
\O_{\ss-style{CY}} \subset \w^n (V\otimes \C)^*\oplus \w^2 V^*.$$

Let $V$ be a real $2n$ dimensional vector space with a Calabi-Yau structure 
$\Phi^0 =(\Ome^0,\ome^0)$. 
We define the complex Hodge star operator $*_\C$ by 
$$
\a \w *_\C \b = <\a \, ,\b> \Ome ^0 ,
$$
where $\a,\b \in \w^{*,0}$. 
The complex Hodge star operator $*_\C$ is a natural generalization of the ordinary Hodge star $*$, 
$$
*_\C \: \w^{i,0} \to \w^{n-i,0}.
$$
The vector space $E^0$ is , by definition, 
$$
E^0_{\ss-style{CY}}(V) = \{\, (i_v \Ome^0, i_v\ome^0 ) \in \w^{n-1,0}\oplus
\w^1\,|
\, v \in V\, \}
$$
The map $TX \to  \w^{n-1,0}$ is given
by $v \mapsto i_v\Ome^0$. Then we see that this map is  an isomorphism.
Hence the projection to the first component defines an isomorphism: 
$$\align
&E^0_{\ss-style{CY}} \arrow  \w^{n-1,0},\\
(i_v\Ome^0,&\, i_v \ome^ 0) \mapsto i_v\Ome^0
\endalign$$
The $E^1_{\ss-style{CY}}$ is the tangent space of Calabi-Yau structures
$\A_{\ss-style{CY}}(X)$.  Hence by (1) and (2) of definition 4-2-1, the vector space
$E^1(V)=E^1_{\ss-style{CY}}(V)$ is the set of 
$ (\a,\b) \in \w^n_\C \oplus \w^2$
satisfying equations 
$$\align
&\a \w \ome^0 + \Ome^0 \w \b =0 , \\
&\a \w \ol{\Ome^0} + \Ome^0 \w \ol{\a} = n c_n \b \w (\ome^0)^{n-1} 
\tag 4
\endalign
$$
Let $P^{p,q}$ be the primitive cohomology group with respect to
$\ome^0$. Then we have the Lefschetz decomposition, 
$$\align
&\a = \a^{\ss-style{n,0}} + \a^{\ss-style{n-1,1}} +\a ^{\ss-style{n-2,0}}\neg\w
\ome^0
\in P^{n,0} \oplus P^{n-1,1}\oplus P^{n-2,0}\negthinspace\w \ome^0 ,\\
&\b = \b^{\ss-style{2,0}} + \b^{\ss-style{1,1}}+ \b^{\ss-style{0,0}}\neg\w
\ome^0 +
\b^{\ss-style{0,2}} \in
 P^{2,0}\oplus P^{1,1}_\R \oplus P^{0,0}\negthinspace\w \ome^0 \oplus 
P^{0,2},
\tag5
\endalign$$
where $\b^{2,0}=\ol{\b^{0,2}}$ and $P^{1,1}_\R$ 
denotes the real primitive forms of type $(1,1)$.
Then equation (4) is written as 
$$\align 
&\a^{\ss-style{n-2,0}}\w \ome \w \ome +\Ome \w \b^{\ss-style{0,2}} =0, 
\tag 6\\
&\a^{\ss-style{n,0}}\w \ol{\Ome} = nc_n \b^{\ss-style{0,0}}\ome^n
\tag7
\endalign$$
Then we see that (6) gives a relation between 
$\a^{n-2,0}$ and $\b^{2,0}$ and (7) also describes a
relation between $\a^{n,0}$ and $\b^{0,0}$. 
Since there is no relation between the primitive parts 
$P^{n-1,1}$ and $P^{1,1}_\R$, the kernel of the projection $E^1_{\ss-style{CY}} \to
\w^{n,0}\oplus
\w^{n-1,1}$ is given by the primitive forms $P^{1,1}_\R$. Hence we have an exact
sequence: 
$$
\CD
0@>>>P^{1,1}_\R @>>>E^1_{\ss-style{CY}}@>>>
\w^{n,0}\oplus\w^{n-1,1}@>>> 0.
\endCD
\tag8$$
The vector space $E^2_{\ss-style{CY}}$ is the subspace of 
$\w^{n,1}\oplus\w^{n-1,2} \oplus \w^3_\R$.  
We also consider the projection to the first component and 
we have an exact sequence: 
$$
0\arrow (\w^{2,1}\oplus\w^{1,2})_\R\arrow E^2_{\ss-style{CY}}
\arrow\w^{n,1}\oplus\w^{n-1,2}\arrow 0,
\tag9$$
where $(\w^{2,1}\oplus\w^{1,2})_\R$ denotes the real part of
$\w^{2,1}\oplus\w^{1,2}$. 
 Let $X$ be a $2n$ dimensional compact K\"ahler manifold. We denote by $\w^{i,j}$ ( global ) differential forms on $X$ of type $(i,j)$.  The real primitive forms of type $(i,j)$ is denoted by $P^{i,j}_\R$. 
Then we have a complex of forms on $X$ by using the exterior derivative $d$: 
$$
\CD 
0@>>> P^{1,1}_\R @>d>> (\w^{2,1}\oplus \w^{1,2})_\R @>d>>\cdots.
\endCD
\tag10
$$
\proclaim{proposition 4-2-3} 
The cohomology groups of the complex (10) are 
respectively given by 
$$
\Bbb P^{1,1}_\R, \quad   (H^{2,1}(X) \oplus H^{1,2}(X))_\R,
$$
where $\Bbb P^{1,1}_\R$ denotes the harmonic and primitive forms.
\endproclaim
\demo{Proof} 
By using K\"ahler identity, we see that 
a closed primitive form of type $(1,1)$ is harmonic.
Hence the first cohomology group of the complex (10) 
is $\Bbb P^{1,1}_\R$. 
Let $q$ be a  real $d$- exact form of type $\w^{(2,1)}\oplus \w^{(1,2)}$.  The
applying
$\pa{\ol{\pa}}$-lemma, we show that $q$ is written as 
$$
q = da,
$$where 
$a = d^* \eta \in \w^{1,1}_\R$ and 
$\eta \in (\w^{2,1}\oplus \w^{1,2})_\R$. We shall show that there exists $k\in \w^1$ 
such that $d^*\eta + dk \in P^{1,1}_\R$. 
By the Lefschetz decomposition,  the three form $\eta$ is written as
$$
\eta = s + \theta\w \ome^0,
$$
where $s\in (P^{2,1}\oplus P^{1,2})_\R, $ and $\theta
\in \w^1_\R$. 
Let $\W$ be the contraction with respect to the K\"ahler form $\ome^0$.
Since $\W$ and $d^*$ commutes, 
$$
\W d^*  \eta = d^* \W \eta = d^* \W ( s+ \theta \sw \ome^0 ) = d^* \theta.
$$
On the other hand,  applying K\"aher identity again, we have 
$$
\W dk = d\W k +\sqrt{-1} d_c^* k = \sqrt{-1}d_c^* k,
$$
where $d_c^* = \pa^* -\ol{\pa}^*$.
Since $k\in \w^1$,  
$$
\align
d_c^* k =& (\pa^* - \ol{\pa}^*) k = \pa^* k^{1,0} 
- \ol{\pa}^*k^{0,1} \\
=&(\pa^* +\ol{\pa}^*) (k^{1,0}- k^{0,1} ).
\endalign
$$
Hence if we define $k$ by 
$$k= \sqrt{-1} ( \theta^{1,0} -\theta^{0,1} ),$$
then 
$$
\W ( d^*\eta + dk ) = d^*\theta + \sqrt{-1} d^* (k^{1,0}
-k^{0,1} )= d^* \theta + ( - d^* \theta^{1,0} 
-d^*\theta^{0,1} )=0.
$$
Hence each exact form $q$ of type $(\w^{2,1}\oplus 
\w^{1,2})_\R$ is given by 
$$
q = d ( d^*\eta + dk),
$$
where $d^* \eta + dk \in P^{1,1}_\R$.
Thus the second cohomology group of the complex (10) 
is $(H^{2,1}(X)\oplus H^{1,2}(X))_\R$.
\qed\enddemo
\proclaim{Theorem 4-2-4} 
The cohomology groups of the complex $\#_{\ss-style{CY}}$:
$$
\CD 
0@>>> E^0_{\ss-style{CY}}@>d>> E^1_{\ss-style{CY}}
@>d>> E^2_{\ss-style{CY}} @>d>>\cdots,
\endCD
$$
is respectively given by 
$$\align
&H^0(\#_{\ss-style{CY}}) = H^{n-1,0}(X), \\
&H^1(\#_{\ss-style{CY}}) = H^{n,0}(X)\oplus H^{n-1,1}(X) \oplus P^{1,1}_\R, \\
&H^2(\#_{\ss-style{CY}}) = 
H^{n,1}(X)\oplus H^{n-1,2}(X) \oplus (H^{2,1}(X)\oplus H^{1,2}(X))_\R,
\endalign
$$
In particular , $p^k$ is injective for $k=0,1,2$. 
\endproclaim
\demo{Proof} 
By (8) and (9),  we have the following diagram: 
$$
\CD 
@.@.0 @.0\\
@.@.@VVV @VVV \\
 @.0 @>>> P^{1,1}_\R @>>> (\w^{2,1}\oplus\w^{1,2})_\R @>>>\cdots \\
@. @VVV @VVV @VVV \\
0@>>>E^0_{\ss-style{CY}}@>>>E^1_{\ss-style{CY}}
@>>>E^2_{\ss-style{CY}}@>>>\cdots \\
@.@VVV @VVV @VVV\\
0@>>>\w^{n-1,0}@>>> \w^{n,0}\oplus\w^{n-1,1}
@>>> \w^{n,1}\oplus \w^{n-1,2}@>>>\cdots \\
@.@VVV @VVV@VVV \\
@. 0 @.0 @.0
\endCD
$$
At first we shall consider H$^2(\#_{\ss-style{CY}})$.
We assume that $(s,t)\in E^2_{\ss-style{CY}}$ is written as an exact form, i.e., $(s,t) = (da, db)$.
Let $a$ be an element of $\w^{n,0}\oplus \w^{n-1,1}$. 
There is a splitting map $\lam\: \w^{n,0}\oplus\w^{n-1,1} \to
\w^2$ such that  
$(a, \lam(a) ) $ is an element of
$E^1_{\ss-style{CY}}$.  Hence 
$$( da , d\lam(a) ) 
\in E^2_{\ss-style{CY}}.
$$
By (10), 
we see that 
$$
db - d\lam(a)\in 
(\w^{2,1}\oplus \w^{1,2})_\R.
$$
Then by proposition 4-2-3,  there exists $p\in P^{1,1}_\R$ such that 
$$
db -d\lam(a) = dp.
$$
Hence $(s,t)$ is written as 
$$
(s,t) =(da ,db ) = ( da, d( \lam(a)+ p ) ) ,
$$
where 
$( a, \lam(a) + p)
\in E^1_{\ss-style{CY}}$. 
Hence we see that 
$$H^1(\#_{\ss-style{CY}}) = H^{n,1}(X)\oplus
H^{n-1,2}(X) \oplus (H^{2,1}(X)\oplus H^{1,2}(X) ) _\R.
$$
Next we shall consider $H^1(\#_{\ss-style{CY}})$. 
Let $(a,b)$ be an element of $E^1_{\ss-style{CY}}$ 
and we assume that $(a,b) = ( d\eta, d\gam )$. 
Then $s$ is written as $s= i_v \Ome^0$ for some 
$v \in TX$. 
By our definition $E^0_{\ss-style{CY}}$, 
$( i_v \Ome^0 ,i_v \ome^0)$ is an element of $E^0_{\ss-style{CY}}$. 
Hence $d\gam - di_v \ome^0 \in P^{1,1}_\R$. 
By proposition 4-2-3,  a $d$-exact, primitive form 
vanishes. Thus 
$dt - di_v \ome^0 =0$. 
Hence $(a,b ) = (d\eta, d\gam ) = ( d i_v\Ome^0, 
di_v \ome^0)$, 
where $( i_v \Ome^0 ,i_v\ome^0) \in E^0_{\ss-style{CY}}$. 
Hence we see that 
$$
H^1(\#_{\ss-style{CY}}) = H^{n,0}(X) \oplus 
H^{n-1,1}(X) \oplus \Bbb P^{1,1}_\R(X).
$$
Similarly we see that $E^0_{\ss-style{CY}}(X) =
H^{n-1,0}(X)$. 
\qed\enddemo
Hence we have 
\proclaim{Theorem 4-2-5} 
The orbit $\O_{\ss-style{CY}}$ is metrical, elliptic and topological.
\endproclaim
We also have
\proclaim{Theorem 4-2-6} 
The cohomology group $H^1(\#)$ is the subspace of 
$H^n(X,\C) \oplus H^2 (X,\R)$ which is defined by equations 
$$\align
&\a \w \ome+ \Ome \w \b =0 , \\
&\a \w \ol{\Ome} + \Ome \w \ol{\a} = n c_n \b \w\ome^{n-1}, 
\endalign
$$
where $\a \in H^n (X,\C) , \b \in H^2 (X,\R)$.
\endproclaim
Let $P^{p,q}(X)$ be the primitive cohomology group with respect to
$\ome$. Then we have Lefschetz decomposition, 
$$
\a = \a^{\ss-style{n,0}} + \a^{\ss-style{n-1,1}} +\a ^{\ss-style{n-2,0}}\neg\w
\ome 
\in P^{n,0}(X) \oplus P^{n-1,0}(X) \oplus P^{n-2,0}(X)\negthinspace\w \ome .
$$
$$
\b = \b^{\ss-style{2,0}} + \b^{\ss-style{1,1}}+ \b^{\ss-style{0,0}}\neg\w \ome +
\b^{\ss-style{0,2}} \in
 P^{2,0}(X)\oplus P^{1,1}(X) \oplus P^{0,0}(X)\negthinspace\w \ome \oplus 
P^{0,2}(X).
$$
Then the equation in theorem 4-2-6 is written as 
$$\align 
&\a^{\ss-style{n-2,0}}\w \ome \w \ome +\Ome \w \b^{\ss-style{0,2}} =0, \\
&\a^{\ss-style{n,0}}\w \ol{\Ome} = nc_n \b^{\ss-style{0,0}}\ome^n
\endalign$$
We see that
$\a^{\ss-style{n,0}} \in P^{n,0}(X)$ and $\b^{\ss-style{0,0}}\in P^{0,0}(X)$ are
corresponding to  the deformation in terms of
constant multiplication: 
$$
\Ome \arrow t \Ome, \quad
\ome \arrow s \ome
$$
If a K\"ahler class $[\ome]$ is not invariant under 
a deformation, such a deformation corresponds to an element of 
$\b^{2,0}$ and $\a^{n-2,0}$.  This is in the case of Calabi family of
hyperK\"ahler manifolds, i.e.,  Twistor space gives such a deformation.
It must be noted that there is no relation between $\a^{n-1,1}\in P^{n-1,1}(X)$ and 
$\b^{1,1}(X)\in P^{1,1}(X)$.
We have from theorem 1-8 in section 1,
\proclaim {Theorem 4-2-7}
The map $P$ is locally injective, 
$$
P\: \M_{\ss-style{CY}}(X) \arrow H^n (X,\C) \oplus H^2(X,\R).
$$
\endproclaim
We also have from theorem 1-9 in section 1,
\proclaim{Theorem 4-2-8} 
Let $I(\Ome,\ome)$ be the isotropy group of $(\Ome,\ome)$, 
$$
I(\Ome,\ome) = \{ \, f \in \text{\rm Diff}_0(X)\, |\, f^*\Ome = \Ome , \, f^*
\ome =
\ome \, \}.
$$
We consider the slice $S_0$ at $\Phi^0=( \Ome^0, \ome^0)$. 
Then the isotropy group $I( \Ome^0,\ome^0)$ is a subgroup of 
$I(\Ome,\ome)$ for each $( \Ome,\ome) \in S_0$.
\endproclaim

We define the map $P_{H^2}$ by 
$$
P_{H^2}\: \M_{\ss-style{CY}}(X) \arrow \Bbb P(H^2 (X)),
$$where
$$
P_{H^2} ( [\Ome,\ome] ) \arrow [\ome]_{dR} \in \Bbb P (H^2 (X)),
$$
$\Bbb P(H^2(X))$ denoted the projective space $(H^2(X)-\{0\})/ \R^*$. 
Then we have 
\proclaim{Theorem 4-2-9} 
The inverse image $P^{-1}_{H^2}([\ome]_{dR})$ 
is a smooth manifold. 
\endproclaim
\demo{Proof}
From theorem 4-2-6 and theorem 4-2-7 the differential of the 
map $P_{H^2}$ is surjective. Hence from the implicit function theorem 
$P^{-1}_{H^2}([\ome]_{\ss-style{dR}})$ is a smooth manifold.
\enddemo
\demo{Remark} 
$P^{-1}_{H^2}([\ome]_{\ss-style{dR}})$ is 
the $\C^*$ bundle over the moduli space of polarized manifolds [8].
\enddemo

\head \S 5 HyperK\"ahler structures
\endhead
Let $V$ be a $4n$ dimensional real vector space. 
A hyperK\"ahler structure on $V$ consists of a metric $g$ and three complex
structures 
$I$,$J$ and $K$ which satisfy the followings: 
$$\align
&g( u, v ) = g( Iu ,Iv) = g( Ju Jv) =g( Ku,Kv),
\quad\text{     for }  u, v \in V,\tag1 \\
&I^2 = J^2 =K^2 = IJK =-1.
\tag2\endalign
$$
Then we have the fundamental two forms $\ome_I, \ome_J, \ome_K$ by 
$$\align 
\ome_I( u,v ) =&g ( Iu, v), \,\,
\ome_J ( u,v) =g(Ju,v), \\
&\ome_K(u,v) = g(Ku,v).
\tag3\endalign $$
We denote by $\ome_\C$ the complex form $\ome_J + \sqrt {-1} \ome_K$. 
 Let $\A_{HK}(V)$ be the set of pairs $(\ome_I,\ome_\C)$ corresponding to 
hyperK\"ahler structures on $V$. As in section one $\A_{HK}(V)$ is the subset of 
$\w^2 \oplus\w^2_\C$ and the group GL$(4n, \R)$ acts on $\A_{HK}(V)$. 
Then we
see that $\A_{HK}(V)$  is GL$(4n ,\R)-$orbit with the isotropy group Sp$(n)$,
$$
\A_{HK}(V) = GL(4n, \R) / Sp(n).
\tag4$$
We denote by $\O_{HK}$ the orbit $\A_{HK}(V)$.
\proclaim{Theorem 5-1}
The orbit $\O_{HK}$ is metrical,elliptic and topological.
\endproclaim
Let $\Phi^0 = ( \ome_I^0 ,\ome_J^0, \ome_K^0 )$ 
be a hyperK\"ahler structure on a $4n$ dimensional vector space $V$.  We denote by $\ome^0_\C$ the complex symplectic form $\ome^0_J + \sqrt{-1}\ome^0_K$.  Then we consider the pair $( \ome^0_I , \ome^0_\C) $.  The vector space $E^k_{\ss-style{HK}}$ are 
respectively given by 
$$\align
&E^0_{\ss-style{HK}}= \{ \, ( i_v \ome^0_I , i_v \ome^0_\C) \, |\, v\in TX \, \}  \\
&E^1_{\ss-style{HK}}= \{ \,  (\hrho_a \ome^0_I, 
\hrho_a \ome^0_\C) \, |\, a\in End (TX) \, \}.
\endalign
$$
Then we consider the projection to the second  component
and we have the diagram: 
$$
\CD
0@>>>E^0_{\ss-style{HK}}@>>>E^1_{\ss-style{HK}}
@>>>E^2_{\ss-style{HK}}@>>>\cdots \\
@.@VVV @VVV @VVV\\
0@>>>\w^{1,0}@>>> \w^{2,0}\oplus\w^{1,1}
@>>> \w^{3,0}\oplus \w^{2,1}\oplus \w^{1,2}@>>>\cdots \\
@.@VVV @VVV@VVV \\
@. 0 @.0 @.0
\endCD
$$
Let $I,J,K$ be the three almost complex structures on $V$.  Then we denote by $\w^{1,1}_I$ forms of type $(1,1)$ with respect to $I$. 
 Similarly $\w^{1,1}_J $( resp. $\w^{1,1}_K$) 
denotes forms of type $(1,1)$ w.r.t $J$ ( resp. $K$).
We define $\w^2_{\ss-style{HK}}$  by the intersection between them, 
$$
\w^2_{\ss-style{HK}}= \w^{1,1}_I\cap \w^{1,1}_J 
\cap \w^{1,1}_K.
$$ 
Note that $a \in \w^2_{\ss-style{HK}}$
is the primitive form with respect to $I,J,$ and $K$. 
When we identify two forms with $so(4m)$, 
we have the decomposition: 
$$
\w^2 = sp(4m ) \oplus so(4m)/ sp(4m).
$$
Then $\w^2_{\ss-style{HK}}$ corresponds to sp$(4m)$. 
Hence the dimension of $\w^2_{\ss-style{HK}}$ is 
$2m^2 + m$. 
We also see that 
$$\align
&\dim_\R E^1_{\ss-style{HK}}= \dim_\R gl(4m,\R) 
/ sp(4m) = 14 m^2 -m, \\
&\dim_\R  \w^{2,0}\oplus \w^{1,1} = 12 m^2 - 2m
\endalign
$$
In fact we see that the kernel of the  map 
$E^1_{\ss-style{HK}} \to \w^{2,0}\oplus \w^{1,1}$ 
is given by $\w^2_{\ss-style{HK}}$ . 
We also define $\w^3_{\ss-style{HK}}$ 
by real forms of type $(\w^{2,1}\oplus \w^{2,1})_\R$
for each $I,J,$ and $K$.
Then we also see that
the  kernel of the map $E^2_{\ss-style{HK}}\to 
\w^{3,0}\oplus\w^{2,1}\oplus\w^{1,2}$
is $\w^3_{\ss-style{HK}}$.
We consider the following complex:
$$
\CD 
0 @>>> \w^2_{\ss-style{HK}}@>>> \w^3_{\ss-style{HK}}@>>> \cdots 
\endCD
\tag HK$$
As in proof of Calabi-Yau structures, 
we see that the cohomology groups of the complex (HK) are respectively  given by 
$$
\align 
&{\Bbb H}^2_{\ss-style{HK}} = 
\{\,\text{real harmonic forms of type} (1,1) w.r.t\, I,J,K \, \}\\
&{\Bbb H}^3_{\ss-style{HK}} =\{ \,
\text{real harmonic forms of type} \w^{2,1}\oplus \w^{1,2}
w.r.t \, I.J.K\, \}.
\endalign
$$
Hence we have the following:
$$
\CD 
@.@.0 @.0\\
@.@.@VVV @VVV \\
 @.0 @>>> \w^2_{\ss-style{HK}} @>>> 
 \w^3_{\ss-style{HK}} @>>>\cdots \\
@. @VVV @VVV @VVV \\
0@>>>E^0_{\ss-style{HK}}@>>>E^1_{\ss-style{HK}}
@>>>E^2_{\ss-style{HL}}@>>>\cdots \\
@.@VVV @VVV @VVV\\
0@>>>\w^{1.0}@>>> \w^{2,0}\oplus\w^{1,1}
@>>> \w^{3.0}\oplus \w^{2,1}\oplus\w^{1,2}@>>>\cdots \\
@.@VVV @VVV@VVV \\
@. 0 @.0 @.0
\endCD
$$
\proclaim{Theorem 5-2} 
The cohomology groups of the complex $\#_{\ss-style{HK}}$  are given by 
$$
\align 
&H^0(\#_{\ss-style{HK}}) = H^{1,0}(X) \\
&H^1(\#_{\ss-style{HK}}) = H^{2,0}(X) \oplus H^{1,1}(X) \oplus \Bbb H^2_{\ss-style{HK}},\\
&H^3(\#_{\ss-style{HK}}) = 
H^{3,0}(X)\oplus H^{2,1}(X)\oplus H^{1,2}(X)\oplus 
\Bbb H^3_{\ss-style{HK}}.
\endalign
$$
In particular, the map $p^k$ is injective for $k=0,1,2$. 
\endproclaim
\demo{Proof} 
The proof is essentially same as in the case of Calabi-Yau structures. 
Let $\lam $ be the splitting map $ \w^{2,0}\oplus \w^{1,1} \to E^1_{\ss-style{HK}}$. 
Let $(s,t)$ by an element of $E^2_{\ss-style{HK}}$. 
We assume that $(s,t) = (da, db)$ for $b \in \w^{2,0}\oplus \w^{1,1}$ and $a \in \w^2$.
By using the splitting map $\lam$, we have 
$(  \lam(b) , b) \in E^1_{\ss-style{HK}}$.  
Then $(  d\lam(b) , db) \in E^2_{\ss-style{HK}}$. 
Hence $da- d\lam (b) \in \w^3_{\ss-style{HK}}$.
Then there is an element $\gam \in \w^2_{\ss-style{HK}}$ such that 
$$
da -d\lam (b) = d\gam .
$$
Hence $(s,t)= ( da, db) =( d(\lam(b) +\gam) , db ) $, 
where
$(( \lam (b) +\gam , b ) \in E^1_{\ss-style{HK}}$.
Thus we have 
$$
H^2(\#_{\ss-style{HK}}) = H^{3,0}(X) \oplus 
H^{2,1}(X)\oplus H^{1,2}(X) \oplus \Bbb H^3_{\ss-style{HK}}.
$$
Similarly we see that 
$$\align
&H^1(\#_{\ss-style{HK}}) = H^{2,0}(X)\oplus H^{1,1}(X)\oplus \Bbb H^2_{\ss-style{HK}}\\
&H^0(\#_{\ss-style{HK}}) = H^{1,0}(X)
\endalign 
$$
\qed\enddemo
\demo{Proof of theorem 5-1} 
This follows from theorem 5-2.
\qed\enddemo

\head \S6. $G_2$ structures
\endhead
Let $V$ be a real $7$ dimensional vector space with a positive definite
metric. We denote by $S$ the spinors on $V$.
Let $\sig^0$ be an element of $S$ with $\| \sig^0 \| =1$. 
By using the natural inclusion $S\otimes S\subset \w^* V^*$,
we have a calibration by a square of spinors, 
$$
\sig^0 \otimes \sig^0 = 1 + \phi^0 + \psi^0 + vol,
$$
where vol denotes the volume form on $V$ and 
$\phi^0$ ( resp. $\psi^0$ ) is called the {\it associative 3 form }
( resp. {\it coassociative 4 form} ).
Our construction of these forms in terms of spinors is 
written in chapter IV \S 10 of [20] and in section 14 of [11]. 
Background materials of $\G$ geometry are found in [13],[15] and
[24]. 
We also have an
another description of
$\phi^0$ and $\psi^0$. We decompose $V$ into a real $6$ dimensional
vector space $W$ and  the  one dimensional vector space $\R$. 
Let $(\Ome^0,\ome^0)$ be an element of Calabi-Yau structure on $W$
and 
$t$ a nonzero $1$ form on $\R$. 
Then the $3$ form $\phi^0$ and the $4$ form $\psi^0$ are respectively
written as 
$$
\phi^0 = \ome^0 \w t +\text{Im }\Ome^0 ,\quad \psi^0 = \frac 12 \ome^0 \w
\ome^0 -\text{Re }\Ome^0 \w t.
$$
Then as in section 1, we define $G_2$ orbit $\O=\O_{\G}$ as 
$$
\O_{G_2} = \{\, ( \phi, \psi ) = ( \rho_g \phi^0, \rho_g \psi^0 )\, |\,
g \in \text{GL}(V) \, \}.
$$
Note that the isotropy group is the exceptional Lie group $\G$.
We denote by $\A_{G_2}(V)$ the orbit $\O_{G_2}$.
Let $X$ be a real $7$ dimensional compact manifold. 
Then we define a $GL(7,\R)/G_2$ bundle $\A_{G_2}(X)$ by 
$$
\A_{G_2}(X) = \underset {x\in X} \to \bigcup \A_{G_2}(T_x X).
$$
Let $\E^1_{G_2}$ be the set of smooth global sections of $\A_{G_2}(X)$, 
$$
\E^1_{G_2} (X) = \Gam ( X, \A_{G_2}(X) ).
$$
Then the moduli space of $G_2$ structures over $X$ is given as
$$
\M_{G_2}(X) = \{ \, ( \phi, \psi ) \in \E^1_{G_2} \, |\, d \phi =0, d \psi =0 
\,\} / \text{Diff}_0(X).
$$
We shall prove unobstructedness of $\G$ structures.
\proclaim{Theorem 6-1}
The orbit $\O_{G_2}$ is metrical, elliptic and topological.
\endproclaim
The rest of this section is devoted to prove theorem 6-1. 
In the case of $G_2$,  each $E^i$ is written as 
$$
\align 
&E^0 = E^0_{G_2} = \{\, (i_v \phi^0 , i_v \psi^0 ) \in \w ^2\oplus \w^3\, |\, v
\in V
\,\}\\ &E^1 =E^1 _{G_2} =\{ \, (\rho_\xi \phi^0 , \rho_\xi \psi^0 ) 
\in \w^3 \oplus \w^4 \, |
\, \xi \in \frak {gl}(V)\, \}\\
&E^2 =E^2_{G_2} = \{ \, (\theta \w \phi,\theta\w \psi )\in \w^4\oplus
\w^5\,|
\, \theta \in \w ^1 , (\phi,\psi) \in E^1_{G_2} \,\}. 
\endalign
$$
The Lie group $G_2 $ is a subgroup of SO$(7)$
and we see that $\G=\{ \, g \in \text{GL}(V) \, |\,
\rho_g \phi^0 =\phi^0 \, \}$. Hence we have the metric
$g_{\phi}$ corresponding to each $3$ form
$\phi$. Let $*_{\phi}$ be the Hodge star operator with respect to the metric 
$g_\phi$. Then a non linear operator $\Theta(\phi)$ is defined as 
$$
\Theta(\phi) = *_\phi \phi.
\tag1$$
According to [13], the differential of $\Theta$ at $\phi$ is described as 
$$
J( \phi) = d\Theta (a)_\phi = \frac43 *\pi_1 (a) + *\pi_7(a)
-*\pi_{27}(a),
\tag2$$ for each $a \in \w^3$,
where we use the irreducible decomposition of $3$ forms on $V$ under the
action of $\G$, 
$$
\w^3 = \w^3_1 +\w^3_7 +\w^3_{27},
\tag3$$
and each $\pi_i$ is the projection to each component for $i=1,7,27$,
( see also [12] for the operator $J$ ).
From (1) the orbit $\O_{\G}$ is written as 
$$
\O_{\G} = \{ \, (\phi, \Theta (\phi) ) \,|\, \phi \in \w^3\, \}.
\tag4$$
Since $E^1_{\G}(V)$ is the tangent space of the orbit $\O_{\G}$ at
$(\phi^0,\psi^0)$, 
from (2) the vector space $E^1_{\G}(V)$ is also written as
$$
E^1_{\G} (V) =\{ \, ( a, Ja ) \in \w^3\oplus \w^4\, |\, a \in \w^3\,\}.
\tag 5$$

 Let $X$ be a real $7$ dimensional compact manifold and 
$(\phi^0,\psi^0)$ a closed element of $\E^1_{\G}(X)$. Then we
have a vector bundle 
$E^i_{\G}(X)\to X$ by 
$$
E^i_{\G}(X) = \underset {x \in X}\to \bigcup E^i_{\G}( T_xX),
\tag6$$
for each $i=0,1,2$. 
Then we have the complex \#$_{G_2}$, 
$$
\CD
0@>>>\Gam ( E^0_{G_2}) @>d_0>>\Gam(E_{\G}^1)@>d_1>>\Gam(E^2_{\G})@>>>\cdots.
\endCD
$$
The complex \#$_{\G}$ is a subcomplex of the de Rham complex, 
$$
\CD 
0@>>>\Gam ( E^0_{G_2}) @>d_0>>\Gam(E_{\G}^1)@>d_1>>\Gam(E^2_{\G})@>d_2>>\cdots\\
@.@VVV @VVV@VVV \\
\cdots@>>>\Gam ( \w^2\oplus \w^3) @>d>>\Gam(\w^3\oplus\w^4)@>d>> 
\Gam (\w^4\oplus \w^5)@>d>>\cdots.
\endCD
$$
Then we have the map $p^1\: H^1 ( \#_{\G} ) \to H^3(X) \oplus H^4 (X)$ 
and $p^2\: H^2(\#_{\G})\to H^4(X)\oplus H^5(X)$.
The following lemma is shown in [12].
\proclaim{Lemma 6-2}
Let $a^3=db^2$ be an exact $3$ form, where $b^2 \in \Gam (\w^2)$.
If $dJdb^2=0$, then there exists $\gam^2 \in \Gam (\w^2_7)$ 
such that $db^2 = d\gam^2$.
\endproclaim
We shall show that $p^1$ is injective by using lemma 6-2.
\proclaim{Proposition 6-3} 
Let $\a =(a^3, a^4)$ be an element of $\Gam (E^1_{\G})$. 
We assume that there exists $(b^2, b^3)\in \Gam (\w^2\oplus \w^3)$
such that 
$$
(a^3, a^4) = (db^2, db^3).
\tag 7
$$
Then there exists $\gam=(\gam^2,\gam^3) \in \Gam(E^0_{\G})$ 
satisfying 
$$
(db^2,db^3) = (d\gam^2, d\gam^3).
$$
\endproclaim
\demo{Proof} 
From $(5)$ an element of $\Gam(E^1_{\G})$ is written as 
$$
(a^3, a^4 ) =( a^3, Ja^3).
$$
From $(7)$ we have 
$$
dJdb^2 =da^4 =ddb^3=0
\tag 8$$

From lemma 6-2 we have $\gam^2 \in \Gam( \w_7^2)$ such that 
$$
db^2 =d\gam^2.
\tag 9$$
Since $\gam \in \Gam (\w_7^2)$, $\gam$ is written as 
$$
\gam = i_v \phi^0,
\tag10$$
where $v $ is a vector field.
Since $\phi^0$ is closed, $d\gam $ is given by the Lie derivative,
$$
d\gam= di_v \phi^0 =L_v \phi^0.
\tag11$$
Then since Diff$_0$ acts on $\E^1_{\G}$, 
$(L_v \phi^0, L_v\psi^0 )= ( di_v \phi^0, di_v \psi^0)$ is an element of 
$\Gam (E^1_{\G})$.
 Hence from $(5)$, we see 
$$
di_v \psi ^0 = J di_v \phi^0 = Jd\gam^2.
\tag 12
$$
From $(12)$ we have 
$$
( db^2, db^3) = (db^2,Jdb^2) = ( di_v\phi, di_v \psi),
\tag13$$
where $(i_v\phi^0, i_v \psi^0 ) \in \Gam(E^0_{\G})$.
\qed\enddemo
Next we shall show that $p^2$ is injective.
\proclaim{Lemma 6-4} Let $V$ be a real $7$ dimensional vector space with 
a $\G$ structure $\Phi^0_V$. 
Let $u$ be a non-zero one form on $V$. Then 
for any two form $\eta$ there exists $\gam \in \w^2_{14}$ such that 
$$\align
&u\w J(u\w\eta) =u\w J(u\w \gam)
= -2*\|u\|\gam,\\
&i_v \gam=0,
\endalign$$
where $v$ is the vector which is metrical dual of the one form $u$ and $*$ is the
Hodge star operator.
\endproclaim
\demo{Proof} 
The two forms $\w^2$ is decomposed into the irreducible representations of $\G$, 
$$
\w^2= \w^2_7 \oplus \w^2_{14}.
$$
We denote by $\eta_{7}$ the $\w^2_{7}$-component 
of $\eta \in \w^2$. The subspace $u{\ss-style{\w}\, } \w^2$ 
is defined by $\{ u \w \eta \in \w^3 | \eta \in  \w ^2\}$. we also denote by $u{\ss-style\w\,} \w^2_7$ the subspace 
$\{ u\w \eta_7 \in \w^3 |\eta\in \w^2\} $.
Then we have the orthogonal decomposition , 
$$
u{\ss-style\w} \w^2 =u{\ss-style\w} \w^2_7 \oplus (u{\ss-style\w} \w^2_7)^\perp, 
\tag6-4-1$$
where $(u{\ss-style\w}\w^2_7)^\perp$ is the orthogonal complement. 
By the decomposition 6-4-1, $u\w \eta$ is written as 
$$
u \w \eta = u\w \eta_7 + u\w \h\eta.
\tag6-4-2$$
for $\h\eta \in \w^2$. 
Then we see that 
$$i_v (u\w \h\eta) \in \w^2_{14}.
\tag6-4-3$$ 
Since $\eta_7$ is expressed as $i_w \phi^0$ for 
$w \in V$, we have 
$$\align
u\w J (u\w \eta_7 ) =& u\w J(u\w i_w \phi^0 ) \\
= &u\w J \hrho_a \phi^0,
\endalign$$
where $a= w \otimes u \in V\otimes V^* \cong 
$End$(V)$.
Since $J\hrho_a \phi^0= \hrho_a\psi^0$,
$$
u\w J\hrho_a \phi^0 = u\w \hrho_a\psi^0 = 
u\w (u\w i_w \psi^0 ) =0.
$$
Hence $$ u\w J ( u\w \eta_7) =0. \tag 6-4-4$$
Then by 6-4-2 we have 
$$
u\w J(u\w \eta) = u\w J ( u\w \h\eta).
\tag6-4-5$$
$\h\eta$ is written as 
$$
\h\eta = \frac 1{ 2 \| u \|^2}( i_v (u \w \h\eta) + 
u\w i_v \h\eta ).
\tag6-4-6$$
We define $\gam$ by 
$$
\gam =\frac1{2\| u\|^2} i_v  (u\w \h\eta).
$$
By 6-4-3, $\gam \in \w^2_{14}$. 
By 6-4-5,6 we have 
$$
u\w J(u\w \eta ) = u\w J( u\w \gam).
\tag6-4-7$$
Since $\gam \in \w^2_{14}$, $\gam \w \psi^0 =0$. 
Then it follows that 
$$
\psi^0 \w u \w \gam =0.
\tag 6-4-8$$
We also have $* \gam =- \gam \w \phi^0$ from 
$\gam \in \w^2_{14}$. 
Since $i_v \gam =0$, we have $u\w (* \gam ) =0$. 
Thus 
$$
 \phi^0 \w u\w \gam  =0.
\tag6-4-9$$
By 6-4-8 and 6-4-9, we have 
$$
u \w \gam \in \w^3_{27}.
\tag6-4-10$$
Then by 6-4-7, 
$$\align
u\w J (u \w \eta ) =& u \w J( u \w \gam ) \\
=& -u\w * ( u\w \gam ) = - * i_v u \w \gam \\
=& -2 \|u \|^2 ( * \gam )
\endalign$$
\qed\enddemo
\proclaim{Proposition 6-5}
 Let $E^2_{\G}(V)$ be the vector space as in before. 
 Then we have an exact sequence, 
$$\CD 
0@>>>\w^5_{14}@>>>E^2_{\G}(V) @>>>\w^4 @>>>0
\endCD
$$
\endproclaim
\demo{Proof} 
The map $E^2_{\G}\to \w^4$ is the projection to 
the first component. We denote by Ker the Kernel of the map $E^2_{\G}\to \w^4$.  We shall show that Ker $\cong \w^5_{14}$.
Let $\{ v_1, v_2,\cdots ,v_7\}$ be an orthonormal basis of $V$. 
We denote by $\{ u^1, u^2, \cdots, u^7\}$ the dual basis of $V^*$.
Let $(s,t)$ be an element of $E^2_{\G}(V)$, where 
$s\in \w^4$ and $t\in \w^5$.
Then we have the following description: 
$$
s= u^1\w a_1 +u^2\w a_2 +\cdots+u^7 \w a_7,
\tag6-5-1$$ 
$$
t=u^1\w J a_1 + u^2\w Ja_2+\cdots+u^7\w Ja_7.
\tag6-5-2$$
where $a_1, a_2,\cdots ,a_7 \in \w^3$ satisfying 
$$
i_{v_l}a_m =0, \forall  l<m.
$$
We assume that $(s,t) \in $Ker. Then $s=0$. 
By 6-5-1, 
we see that $u^l \w a_l =0$, for all $l$. 
Hence each $a_l$ is written as 
$$a_l = u^l \w \eta_l
\tag6-5-3$$
where $\eta_l \in \w^2$. 
By (6-5-2) we have 
$$
t=\sum_{l=1}^7 u^l \w J (u^l\w \eta_l).
$$
Then it follows from lemma 6-4 there exists $\gam_l$ such that
$$
t = \sum _{l=1}^7 u^l\w J (u^l\w \gam_l) = 
-2\sum_{l=1}^7  \| u^l \|^2 ( * \gam _l ),
$$
where $\gam_l \in \w^2_{14}$.
Hence $t \in \w^5_{14}$. 
Therefore we see that Ker = $\w^5_{14}$.
\qed\enddemo
\proclaim{Lemma 6-6} 
Let $X$ be a compact $7$ dimensional manifold with 
$\G$ structure $\Phi^0$, (i.e., $d\Phi^0=0).$
Then for any two form $\eta$ there  exists 
$\gam \in \w^2_{14} $ such that 
$$\align
&dJ d\eta = dJ d \gam = -* \trian \gam, \\
&d^* \gam =0.
\endalign $$
\endproclaim
\demo{Proof} 
We denote by $d\w^2$ the closed subspace 
$\{ d \eta | \eta \in \w^2 \}$. 
Since $d\w^2_7 = \{ d \eta_7 |\eta_7 \in \w^2_7 \}$ is the closed subspace of $d\w^2$, we have the decomposition, 
$$
d\w^2 = d\w^2_7  \,\oplus\,\, ( d\w^2_7)^\perp
\tag6-6-1$$
where $( d\w^2_7)^\perp $ denotes the orthogonal subspace of $d\w^2_7$. 
By 6-6-1, $d\eta$ is written as 
$$
d\eta = d\eta_7 + d\h\eta,
$$
where $d\h\eta \in ( d\w^2_7)^\perp$.
Hence we have 
$$
d^* d\h\eta \in \w^2_{14}
\tag6-6-2$$
As in the proof of lemma 6-4, $\eta_7 $ is written as 
$i_w \phi^0$ for some $w \in TX$. Hence 
$$
dJ d \eta_7 = dJ di_w \phi^0 = dJ L_w \phi^0 = 
d L_w \psi^0 = d d i_w \psi^0 =0.
\tag6-6-3$$
Thus $dJ d\eta = dJ d\h\eta$.
By the Hodge decomposition, 
we have 
$$
\h\eta = Harm( \h\eta) + d d^* G\h\eta 
+ d^*d G\h\eta,
$$
where $Harm(\h\eta)$ is the harmonic part of $\h\eta$ and $G$ denotes the Green operator. 
We define $\gam$ by 
$$
\gam = d^* d G \h\eta.
$$Then by Chern's theorem ( $\pi_7 G = G \pi_7 )$ and 6-6-2, we see that $\gam \in \w^2_{14}$.
Then $d\h\eta = d\gam$ and $d^* \gam =0$.
Since $\gam \in \w^2_{14}$, we have 
$\gam \w \psi^0 =0$ and $*\gam = -\gam \w\phi^0$. Hence we have 
$$\align
&d\gam \w \phi^0 =0, \tag6-6-4\\
&d\gam \w \psi^0 =0.\tag6-6-5
\endalign$$
Hence it follows from 6-6-4,5 that 
$$
d\gam \in \w^3_{27}.
\tag6-6-6
$$ 
Then by 6-6-6,
$$
dJd\gam = -d * d\gam = - *\trian \gam.
$$
By 6-6-3, 
$$
dJd\eta = -*\trian \gam .
$$
\qed\enddemo
\proclaim{Proposition 6-7} 
$$
H^2(\#_{\G} ) = H^4 (X) \oplus H^5_{14}(X) .
$$
In particular,  
$$
p^2\: H^2( \#_{\G} ) \arrow H^4 (X)\oplus H^5(X)
$$
is injective.
\endproclaim 
\demo{Proof} 
Let $(s,t)$ be an element of $E^2_{\G}(X)$. 
We assume that $s,t$ are exact forms respectively,i.e., 
$$
s= da, \,t = db,
\tag6-7-1$$
for some $a\in \w^3$ and $b \in \w^4$. 
Then we shall show that there exists $\til{a}\in\w ^3$
such that 
$ da = d\til{a} $ and $db = d J \til{a}$.
Since $( da, dJa)$ is an element of $E^2_{\G}$, 
it follows from proposition 6-5 that 
$$
db - dJa \in \w^5_{14}. \tag 6-7-2
$$ 
We shall show that there exists $\eta \in \w^2$ satisfying, 
$$
db = d  J ( a+ d \eta)
\tag6-7-3$$
In order to solve the equation (6-7-3), we apply lemma 6-6. 
Then there exists $\gam \in \w^2_{14}$ such that 
$$\align
&dJ d\eta = -*\trian \gam \tag 6-7-4\\
&d^* \gam =0.
\endalign
$$
Substituting 6-7-4 to the equation (6-7-3), we have 
$$
-* \trian \gam = db - dJ a
\tag 6-7-5
$$
Then by (6-7-2), there exists a solution $\gam$ of the equation (6-7-5), 
$$
\gam = - G * ( db - dJ a ) \in \w^2_{14}.
$$
Hence  if we set $\til{a} = a + d\gam$, $(s,t)$ is written as 
$$\align
&s = d\til {a} = d( a+ d\gam), \\
&t = dJ\til{a} = dJ ( a+ d\gam )
\endalign$$
Therefore $p^2\: H^2(\#_{\G}) \to H^4(X)\oplus H^5 (X)$ is injective. 
Furthermore we consider harmonic forms 
$\H^4(X)$ and $\H^5_{14}(X)$. 
By Chern's theorem 
$H^4(X) \oplus H^5_{14}(X) \cong \H^4(X)\oplus\H^5_{14}(X)$. 
Since the complex $\#_{\G}$ is elliptic, H$^2(\#_{\G})$ is represented by harmonic forms of the complex $\#_{\G}$,.i.e., 
$$
H^2(\#_{\G})\cong \H^2(\#_{\G})
$$
Then we see that there is the injective map 
$$
\H^4(X) \oplus \H^5_{14}(X) \to \H^2(\#_{\G}). 
$$
Since $p^2$ is injective, we have 
$$
H^2(\#_{\G}) \cong H^4(X) \oplus H^5_{14}(X).
$$
\qed\enddemo
\demo{proof of theorem 6-1} 
By proposition 6-3 and proposition 6-7, 
we have 
$$\align
&H^0(\#_{\G}) \cong H^2_7(X)\cong H^3_7(X) \\
&H^1(\#_{\G}) \cong H^3(X) \cong H^4(X) \\
&H^2(\#_{\G}) \cong H^4(X)\oplus H^5_{14}(X).
\endalign$$
Hence we have the result.
\enddemo

\head \S 7.  Spin$(7)$ structures
\endhead Let $V$ be a real $8$ dimensional vector space with a positive definite metric. 
We denote by $S$ the spinors of $V$. Then $S$ is decomposed into the positive spinor
$S^+$ and the negative spinor $S^-$.  Let $\sig^+_0$ be a positive spinor with $\|\sig_0^+
\| =1$.  Then under the identification $S\otimes S \cong \w^*V$,  we have a calibration by
the square of the spinor, 
$$
\sig^+_0\otimes\sig^+_0 = 1 + \Phi^0 + vol,
$$ where vol denotes the volume form on $V$ and 
$\Phi^0$ is called {\it the Cayley 4 form } on $V$   (see [11], [20] for
our construction in terms of spinors). Background materials of
Spin$(7)$ geometry  are found in [14],[15] and [24]. we decompose
$V$ into  a real $7$ dimensional vector space $W$ and the one
dimensional vector space $\R$, 
$$ V= W \oplus \R.
$$ Then a Cayley $4$ form $\Phi^0$ is defined as 
$$
\phi^0 \w \theta + \psi^0  \in \w^4 V^* ,$$ where $(\phi^0,\psi^0) \in \O_{\G}(W)$
and $\theta$ is non zero one form on
$\R$.  We define an orbit $\O_{Spin(7)}=\A_{\Spin}(V)$ by 
$$
\O_{\Spin} = \{\, \rho_g \Phi^0 \, |\, g \in \text{GL}(V) \, \}.
$$ Since the isotropy is $\Spin$, the orbit $\O_{\Spin}$ is written as 
$$
\O_{\Spin} = GL(V)/\Spin.
$$ Let $X$ be a real $8$ dimensional compact manifold. Then we define $\A_{\Spin}(X)$
by 
$$
\A_{\Spin}(X) = \underset{x\in X}\to \bigcup 
\A_{\Spin}(T_x X) \arrow X.
$$ We denote by $\E^1_{\Spin}(X)$ the set of global section of $\A_{\Spin}(X)$, 
$$
\E^1_{\Spin}(X) = \Gam (X, \A_{\Spin}(X) ).
$$ Then we define the moduli space of $\Spin$ structures over $X$ as 
$$
\M_{\Spin}(X) = \{ \, \Phi \in \E^1_{\Spin}\, |\, d \Phi =0 \, \}/
\text{Diff}_0(X).
$$ The following theorem is shown in [15]
\proclaim{Theorem 7-1} [15] The moduli space $\M_{\Spin}(X)$
is a smooth manifold with 
$$
\dim \M_{\Spin}(X) = b^4_1 +b^4_7+ b^4_{35},
$$ where Harmonic $4$ forms on $X$ is decomposed into irreducible representations of
Spin $(7)$,
$$
\H^4 (X) = \H^4_1\oplus \H^4_7 \oplus\H^4_{27} \oplus\H^4_{35},
$$ each $b^4_i$ denoted $\dim \H^4_i$, for $i = 1,7,27$ and $35$.
\endproclaim Note that $\H^4(X)$ is decomposed into self dual forms and anti-self dual
forms, 
$$
\H^4 (X) = \H^+\oplus \H^-,
$$ where 
$$
\H^+ (X) = \H^4_1\oplus \H^4_7 \oplus\H^4_{27}, 
\quad \H^- = \H^4_{35}.
$$ We shall show theorem 7-1 by using our method in section
one. 
\proclaim{Theorem 7-2} The orbit $\O_{Spin(7)}$ is metrical, elliptic
and topological.
\endproclaim Since $\Spin$ is a subgroup of SO$(8)$, we have the metric $g_{\phi^0}$ 
for each $\Phi^0 \in \O_{\Spin}$. For each $\Phi^0 \in \O_{\Spin}(V)$, $\w^3$ and
$\w^4$ are orthogonally decomposed into the irreducible representations of $\Spin$, 
$$
\w^3 = \w^3_8 \oplus \w^3_{48},
$$
$$
\w^4 = \w^+ \oplus \w^- = (\w^4_1 \oplus \w^4_7 \oplus \w^4_{27}) 
\oplus \w^4_{35},
$$ where $\w^p_i$ denotes the irreducible representation of $\Spin$  of $i$ dimensional.
We denote by $\pi_i $ the orthogonal projection to each component. Let $X$ be a real $8$
dimensional compact manifold with a closed form
$\Phi^0 \in 
\E^1_{\Spin}(X)$.  Let $g_{\Phi^0}$ be the metric corresponding to $\Phi^0$. Then
there is a unique parallel positive spinor $\sig^+_0
\in \Gam ( S^+)$ with 
$$
\sig^+_0 \otimes \sig^+_0 = 1 + \Phi^0 + \text{vol},
$$ where $S^+ \otimes S^+$ is identified with the subset of Clifford algebra Cliff
$\cong \w^*$ ( see [16]).  By using the parallel spinor $\sig^+_0$, the positive and
negative spinors are respectively identified with following representations,
$$\align
\Gam ( S^+ ) &\cong \Gam ( \w^4_1\oplus \w^4_7 ),\\
\sig^+& \arrow \sig^+ \otimes \sig^+_0,
\tag 1\endalign$$
$$\align
\Gam (S^-) &\cong \Gam ( \w^3_8),\\
\sig^- &\arrow \sig^-\otimes \sig^+_0,
\tag 2\endalign$$ where $\sig^\pm \in \Gam (S^\pm)$. Under the identification (1) and
(2),  The Dirac operator $D^+ \:\Gam (S^+ )\to \Gam (S^-)$ is written as 
$$
\pi_8 \circ d^* \:\Gam ( \w^4_1 \oplus \w^4_7 ) \arrow \Gam (\w^3_8).
$$ In particular Ker $\pi_8 \circ d^*$ are Harmonic forms in $\Gam (\w^4_1\oplus
\w^4_7)$.  Hence we have 
\proclaim{Lemma 7-3} 
$$
\text{\rm Ker }\pi_8\circ d^* = \H^4_1(X) \oplus \H^4_7 (X).
$$
\endproclaim
In the case of $\Spin$, each $E^i
=E^i_{\Spin}$ is given by 
$$ E^0_{\Spin} = \w^3_8, \quad E^1_{\Spin} = \w^4_1\oplus \w^4_7\oplus \w^4_-.
$$ Let $\a$ be an element of $\Gam (E^1_{\Spin}(X))$.  We assume that 
$$ d\a =0 ,\quad \pi_8 d^* \a =0,
\tag3$$ So that is, $\a$ is an element of $\H^1(\#)$, where \# is the complex 
$$
\CD  0@>>>\Gam ( E^0_{\Spin}) @>d_0>>\Gam (E^1_{\Spin})@>d_1>>\Gam
(E^2_{\Spin}) @>>>\cdots\\ @.@|@|@| \\
\cdots@>>>\Gam (\w^3_8 ) @>d>>\Gam (\w^4_1\oplus \w^4_7 \oplus
\w^-)@>d>>\Gam(\w^5)@>>>\cdots
\endCD
$$ (Note that $d_0^* =  \pi_8 d^*$.) We decompose $\a$ into the self-dual form and the
anti-self-dual form, 
$$
\a = \a^+  + \a^-\in \Gam (\w^+) \oplus \Gam (\w^-).
$$ From (3) we have 
$$\align &d\a^+  + da^- =0\\ &\pi_8 *d \a^+ -\pi_8 * d\a^- =0.
\endalign
$$ Hence we have $\pi_8 d^* \a^+ =0$. From lemma 7-3, we see
that $d\a^+=0$.  Hence we also have $d\a^-=0$ and it implies that
$\a$ is a harmonic form with respect to the metric $g_{\Phi^0}$.
Hence the map 
$p\: H^1( \#) \cong \H^1 (\#) \arrow H^4 (X) \cong \H^4(X)$  is
injective.
\proclaim{Theorem 7-4}  The cohomology groups of the complex
$\#_{{\ss-style{Spin(7)}}}$  are  respectively given by 
$$ \align
&H^0(\#_{{\ss-style{Spin(7)}}}) \cong H^3_8(X), \\
&H^1(\#_{{\ss-style{Spin(7)}}}) \cong H^4_1(X)\oplus H^4_7(X)
\oplus H^4_-(X),\\
&H^2(\#_{{\ss-style{Spin(7)}}})= H^5(X),
\endalign
$$
In particular $p^1$ and $p^2$ are respectively injective.
\endproclaim
\demo{Proof} 
It is sufficient to show that $H^2(\#_{{\ss-style{Spin(7)}}})= H^5(X)$.
Since anti-self dual forms $\w^4_-$ is the subset of 
$E^1_{{\ss-style{Spin(7)}}}$, we see that our result.
\enddemo
\demo{Proof of theorem 7-2} 
This follows from theorem 7-4.
\enddemo


\enddocument